\documentclass{amsart}
\let\<=\langle
\let\>=\rangle
\let\noi=\noindent
\let\sse=\subseteq
\let\vphi=\varphi
\let\veps=\varepsilon
\let\limply=\Longrightarrow
\let\what=\widehat

\def\0{\{0\}}
\def\diag{{\rm diag}}
\def\span{{\kern.5pt{\rm span}\,}}
\def\hotimes{{\kern2pt\what\otimes\kern1pt }}
\def\smallfrac#1#2{{\textstyle{\frac{#1}{#2}}}}

\def\conv{{\;\longrightarrow\;}}
\def\wconv{{{\buildrel_{\scriptstyle w}\over\conv}}}
\def\sconv{{{\buildrel_{\scriptstyle s}\over\conv}}}
\def\uconv{{{\buildrel_{\scriptstyle u}\over\conv}}}
\font\fiverm=cmr5
\def\sslash{\hbox{{\fiverm /}}}
\def\notuconv{{{\uconv\kern-13pt\sslash}\kern9pt}}
\def\notsconv{{{\sconv\kern-13pt\sslash}\kern9pt}}
\def\notwconv{{{\wconv\kern-13pt\sslash}\kern9pt}}

\def\A{{\mathcal A}}
\def\B{{\mathcal B}}
\def\C{{\kern1pt\mathcal C}}

\def\H{{\mathcal H}}
	\def\Ie{{\mathcal I}}
\def\K{{\mathcal K}}
	\def\Le{{\mathcal L}}
\def\M{{\mathcal M}}
\def\N{{\mathcal N}}
	\def\Oe{{\mathcal O}}
\def\R{{\mathcal R}}
	\def\Se{{\mathcal S}}
\def\T{{\mathcal T}}
\def\U{{\mathcal U\kern.5pt}}
\def\W{{\mathcal W}}
\def\X{{\mathcal X}}
\def\Y{{\mathcal Y}}
\def\Z{{\mathcal Z}}
\def\BX{{\B[\X]}}
\def\BH{{\B[\H]}}
\def\BK{{\B[\K]}}
\def\BL{{\B[\Le]}}
\def\LX{{\Le[\X]}}

\def\CC{{\mathbb C\kern.5pt}}
\def\FF{{\mathbb F}}
\def\DD{{\mathbb D\kern.5pt}}
\def\JJ{{\mathbb J\kern.5pt}}
\def\QQ{{\mathbb Q\kern.5pt}}
\def\NN{{\mathbb N\kern.5pt}}
\def\RR{{\mathbb R\kern.5pt}}
\def\TT{{\mathbb T\kern.5pt}}
\def\ZZ{{\mathbb Z\kern.5pt}}

\newsymbol\varnothing 203F
\let\void=\varnothing

\def\smallmatrix#1{\null\,\vcenter{
                   \baselineskip=8pt\mathsurround=0pt\ialign{
                   \hfil ${\scriptstyle##}$
                   \hfil &&
                   \hfil ${\scriptstyle##}$
                   \hfil \crcr
                   \mathstrut \crcr
                   \noalign{\kern-\baselineskip}#1 \crcr
                   \mathstrut \crcr
                   \noalign{\kern-\baselineskip} \crcr }}\!}

\theoremstyle{definition}

\begin{document}

\vglue-68pt\noi
\hfill{\it Rendiconti del Circolo Matematico di Palermo}\/,
{\bf 75}(5) (2026) 1--37

\vglue20pt
\title[Weak Stability of Operators]
      {An Exposition on Weak Stability of Operators}
\author{C.S. Kubrusly}
\address{Catholic University of Rio de Janeiro, Brasil}
\email{carlos@ele.puc-rio.br}
\renewcommand{\keywordsname}{Keywords}
\date{January 25, 2026}

%%%%%%%%%%%%%%%%%%%%%%%%%%%%%%%%%%%%%%%%%%%%%%%%%%%%%%%%%  ABSTRACT
\begin{abstract}
This is an expository survey on weak stability of bounded linear
operators acting on normed spaces in general and, in particular, on Hilbert
spaces.\ The paper gives a comprehensive account of the problem of weak
operator stability, containing a few new results and some unanswered
questions.\ It also gives an updated review of the literature on the weak
stability of operators over the past sixty years, including present-day
research trends.\ It is verified that the majority of the weak stability
literature is concentrated on Hilbert-space operators.\ We discuss why this
preference occurs and also why the weak stability of unitary operators
is central to the Hilbert-space stability~problem.
\end{abstract}

\maketitle

\vglue-10pt\noi
%%%%%%%%%%%%%%%%%%%%%%%%%%%%%%%%%%%%%%%%%%%%%%%%%%%%%%%%%  SECTION 1
\section{Introduction}

If one types now ``weak stability'' or ``weakly stable'' in MathSciNet, under
the search term `Title' only, there will appear about 180 entries.\ The term
`stability' (weak or not) applies to a variety of instances in mathematics,
ranging from mathematical logic to operator theory, through ordinary and
partial differential equations, difference and integral equations, dynamical
systems and ergodic theory, numerical analysis and mathematical programming,
system and control theory, and beyond.\ In this paper, we stick to operator
theory --- more specifically, to single operators.\ Thus only part of the
above-mentioned 180 entries refers to the problem considered here, namely,
the weak stability of single operators as posed in Section 3 below.\

\vskip5pt
This is an expository survey on weak stability of bounded linear operators
acting on normed spaces and, especially, on Hilbert spaces$.$ It is, on the
one hand, an exposition that brings weak stability of operators from first
principles to recent research topics and some open questions, and, on the
other hand, it surveys the current literature dealing with several aspects of
weak stability of operators, where the available printed material is mostly
(but not all) on Hilbert-space operators.\

\vskip5pt
There are many reasons for considering Hilbert-space operators in this
context$.$ Perhaps one of the most important features of Hilbert space that
motivates such a preference is complementation --- {\it among Banach spaces,
only Hilbert spaces, up to topological isomorphisms, are complemented}\/ (see
\cite{LT,Kal}$\kern.5pt).$ Then {\it the notion of reducing subspaces}\/ (as
defined in Section 2 below) {\it is exclusive to Hilbert spaces}\/, and this
concept is a significant tool for analysing the stability of operators$.\!$ We
carry on this exposition by starting with normed-space operators, going to
Banach spaces whenever completeness is required, and going into Hilbert space
to cope with many techniques depending on operator decomposition based on
reducing subspaces.\

\vskip5pt
The paper is organised into eight more sections$.$ Section 2 gives a detailed
account of weak convergence, mainly for normed-space operators$.$ Sections 3
and 4 introduce notation and terminology and pose the operator stability
problem, with emphasis on the weak stability of bounded linear operators
acting on Banach and Hilbert spaces, which is the central topic of the
paper$.$ Sections 5, 6, 7, 8, and 9 consist of a review of the weak stability
problem of operators, following the setup built in Section 4$.$ These
sections are listed below.\

\vskip5pt\noi
Section 2. Weak Convergence.\
\vskip2pt\noi
Section 3. Basic Notation and Terminology.\
\vskip2pt\noi
Section 4. Notation and Terminology of Operator Stability.\
\vskip2pt\noi
Section 5. General Aspects of Weakly Stable Operators.\
\vskip2pt\noi
Section 6. Ces\`aro Means and Weak Stability.\
\vskip2pt\noi
Section 7. Weak Stability for a Class of $\,2\!\times\!2\,$ Operator
Matrices.\
\vskip2pt\noi
Section 8. Boundedly Spaced Subsequences and Weak Quasistability.\
\vskip2pt\noi
Section 9. Weak Stability of Unitary Operators.\

\vskip5pt
Section 5 focuses on general results about weak stability, having as a
starting point the Foguel decomposition of Hilbert-space contractions, while
the remaining sections will focus on particular techniques or on particular
classes of operators.\

\vskip5pt
Section 6 deals with applications of Ces\`aro means techniques to weak
stability.\

\vskip5pt
Section 7 considers weak stability and weak quasistability for two cases of
Hilbert-space operators decomposed into a two-by-two operator matrix acting on
an invari\-ant-subspace decomposition$.$ The section is split into three
subsections: Subsection 7.1 deals with Brownian-type operators; Subsection 7.2
presents the notion of weak quasistability; and Subsection 7.3 is concerned
with the Foguel operator.\

\vskip5pt
Section 8 considers the concept of boundedly spaced subsequences of power
sequences of operators and applies it to the notion of weak quasistability,
where the relationship between weak stability and weak quasistability is
investigated.\

\vskip5pt
Section 9 closes the paper by reviewing the current characterisations of weak
stability for Hilbert-space unitary operators.\

\vskip5pt
Since Sections 5 through 9 intend to review the literature on weakly stable
operators, references to propositions or corollaries will be followed by the
publication year$.$ Some results in those   sections are either new, or we
could not find any mention of them in literature; some might even go without
proof, some might not, and in such cases (and only in such cases), statements
are accompanied by proof.\

\vskip5pt
On the other hand, the normed-space results of Section 2 are all accompanied
by a proof$.$ In fact, we could not trace them back to literature$.$ Their
counterparts found in the literature are developed for operators acting on
Hilbert spaces under the umbrella of complementation, orthogonality, reducing
subspaces, or Riesz representation and Fourier series theorems$.$ Thus
Proposition 2.6, regarding weak stability for infinite-dimensional
normed-space operators, seems to be introduced here in such a general form$.$
Naturally, it is chronologically posterior to those results reviewed in
Section 5 and extends some of them to normed-space operators.\

%%%%%%%%%%%%%%%%%%%%%%%%%%%%%%%%%%%%%%%%%%%%%%%%%%%%%%%%%  SECTION 2
\section{Weak Convergence}

Let $\X$ and $\Y$ be normed spaces over a scalar field $\FF$ which, in this
paper, is either the real field $\RR$ or the complex field $\CC.$ Let
$\Le[\X,\Y]$ denote the linear space of all linear transformations of $\X$
into $\Y$, and let $\B[\X,\Y]$ denote the normed space of all bounded linear
transformations of $\X$ into $\Y.$ Set ${\LX=\Le[\X,\X]}$, the algebra of all
linear transformations of $\X$ into itself$.$ By an operator on a normed space
$\X$ (or a normed-space operator) we mean a bounded linear transformation of
$\X$ into itself$.$ Set ${\BX=\B[\X,\X]}$, the normed algebra of all operators
on $\X.$ For any normed space $\X\kern-1pt$, let the Banach space
${\X^*\!=\B[\X,\FF]}$ stand for the dual of $\X\kern-1pt.$ The kernel and
range of ${T\kern-1pt\in\LX}$ will be denoted by $\N(T)$ and $\R(T)$,
respectively.

%%%%%%%%%%%%%%%%%%%%%%%%%%%  DEFINITION 2.1
\vskip5pt\noi
{\bf Definition 2.1.}\
Let $\X$ be a normed space$.$ An $\X$-valued sequence $\{x_n\}$ {\it converges
weakly}\/ (or $\{x_n\}$ is {\it weakly convergent}\/) if the scalar-valued
sequence $\{f(x_n)\}$ converges for every ${f\in\X^*}\!$ in the following
sense$.$ There exists ${z\in\X}$ for which the scalar-valued sequence
$\{f(x_n)\}$ converges to $f(z)$ for every ${f\in\X^*}\!$.

\vskip5.5pt
The vector ${z\in\X}$ is the {\it weak limit}\/ of $\{x_n\}$, which is
unique$.$ In fact, if ${z_1,z_2\in\X}$ are weak limits of $\{x_n\}$, then
${f(z_1-z_2)=0}$ for all ${f\in\X^*}\!$, so that ${z_1-z_2=0}$.

\vskip5.5pt
{\it Warning}\/.\
The above definition of weak convergence in a normed space $\X$ is quite
standard (see, e.g., \cite[Definition 14.1]{BN} or
\cite[Definition 1.13.2]{Meg}$\kern.5pt$).\ Weak convergence will be denoted
by ${x_n\!\wconv z}$, or by ${z=w\hbox{\,-}\lim_nx_n}$, meaning that there
exists ${z\in\X}$ for which ${f(x_n)\to f(z)}$ for every ${f\in\X^*}\!.$ It is
worth noting that ``every normed space has a topology $\T\!$ [the weak
topology on it] such that a sequence in the space converges weakly to an
element of the space if and only if the sequence converges to that element
with respect to $\T.$ For the moment, the statement that a sequence converges
weakly to a certain limit should not be taken to imply anything more than is
stated in Definition 2.1.'' \cite[p.116]{Meg}$.$ (We shall not deal with weak
topology techniques in this section, and so Definition 2.1 may be thought of,
for the moment, as a nontopological definition).

%%%%%%%%%%%%%%%%%%%%%%%%%%%  REMARK 2.2
\vskip5pt\noi
{\bf Remark 2.2.}\
\qquad
${x_n\!\wconv z\;\limply\;\sup_n\|x_n\|\kern-1pt<\kern-1pt\infty}$.

\vskip5pt\noi
Indeed, as is well known, (i) $\X$ is isometrically embedded in its second
dual $\X^{**}$ via the natural mapping ${\varPhi\!:\X\to\X^{**}}$ that assigns
to each vector ${x\in\X}$ the functional
${\varphi_x\in\X^{**}\!=\B[\X^*\!,\FF]}$ given by ${\varphi_x(f)=f(x)}$ for
every ${f\in\X^*}\!$, so that ${\|\varphi_x\|=\|x\|}.$

\vskip4pt\noi
(ii) Thus ${|\varphi_{x_n}(f)|\kern-1pt=\kern-1pt|f(x_n)|}$ for each
nonnegative integer $n$, for every ${f\kern-1pt\in\X^*\!}$, where every
scalar-valued sequence $\{f(x_n)\}$ converges (because ${x_n\!\wconv z}$),
so that every $\{f(x_n)\}$ is bounded, and hence
${\sup_n|\varphi_{x_n}(f)|\kern-1pt=\sup_n|f(x_n)|\kern-1pt<\kern-1pt\infty}$
for \hbox{every ${f\in\X^*\!}.$}

\vskip4pt\noi
(iii) So (as each ${\varphi_{x_n}\!\in\B[\X^*\!,\FF]}$ and $\X^*\!$ is a
Banach space), the Banach--Steinhaus Theorem ensures that
${\sup_n\|\varphi_{x_n}\|\!<\!\infty}$, and therefore
${\sup_n\|x_n\|=\sup_n\|\varphi_{x_n}\|\kern-1pt<\kern-1pt\infty}.$

\vskip5.5pt
Now let $\{T_n\}$ be an $\LX$-valued sequence$.$

%%%%%%%%%%%%%%%%%%%%%%%%%%%  DEFINITION 2.3
\vskip5pt\noi
{\bf Definition 2.3.}\
$\kern-1pt$We say that $\kern-1pt\{T_n\}\kern-1pt$ {\it converges weakly}\/
(or $\kern-1pt\{T_n\}\kern-1pt$ is {\it weakly convergent}\/) if for every
${x\kern-1pt\in\kern-1pt\X}$ the $\X$-valued sequence
$\kern-1pt\{T_nx\}\kern-1pt$ converges weakly in the sense of
Definition 2.1$.\kern-1pt$ That is, for each ${x\kern-1pt\in\kern-1pt\X}$
there is a ${z_x\kern-1pt\in\kern-1pt\X}$ such that
${T_nx\kern-1pt\wconv\kern-.5ptz_x}$ (equivalently, for each
${x\kern-1pt\in\kern-1pt\X}$ there exists a ${z_x\kern-1pt\in\kern-1pt\X}$ for
which ${f(T_nx)\to\kern-1ptf(z_x)}$
\hbox{for every ${f\kern-1pt\in\kern-1pt\X^*}\kern-1pt).$}

%%%%%%%%%%%%%%%%%%%%%%%%%%%  PROPOSITION 2.4
\vskip5pt\noi
{\bf Proposition 2.4.}\
{\it Let\/ $\X$ be a normed space, and let\/ $\{T_n\}$ be a sequence of linear
transformations of\/ $\X$ into itself$.$ The following assertions are
equivalent}\/.
\begin{description}
\item{$\kern-6pt$\rm(a)}
$\{T_n\}$ {\it converges weakly}\/.
\vskip3pt
\item{$\kern-6pt$\rm(b)}
{\it There exists\/ ${T\in\LX}$ such that\/ ${T_nx\wconv Tx}$ for every}\/
${x\in\X}.$
\end{description}

\vskip-2pt\noi
\proof
Consider Definition 2.3$.$ By the uniqueness of the weak limit (as in
Definition 2.1), for each ${x\in\X}$ there is a unique ${z_x\in\X}$ such that
${T_nx\wconv z_x}.$ This ensures the existence of a transformation
${T\!:\X\to\X}$ given by ${Tx=z_x=w\hbox{\,-}\lim_nT_nx}$ for every
${x\in\X}$, such that ${T_nx\wconv Tx}$ every ${x\in\X}.$ That $T$ is linear
comes from the linearity of each $T_n$ and ${f\in\X^*}$, and also from the
linearity of the limiting operation.                                     \qed

\vskip5.5pt
Assertion (b) in Proposition 2.4 says that there exists ${T\kern-1pt\in\LX}$
for which
\goodbreak\vskip4pt\noi
$$
T_nx\wconv Tx
\quad\hbox{for every}\quad
x\in\X.
$$
\vskip1.5pt\noi
According to Definition 2.1, this means that there exists ${T\kern-1pt\in\LX}$
such that
$$
\hbox{for each}\quad
x\in\X,
\quad
f(T_nx)\to f(Tx)
\quad\hbox{for every}\quad
f\in\X^*\!
$$
or, equivalently (since $f$ is linear),
$$
f\big((T_n-T)x\big)\to0
\quad\hbox{for every}\quad
x\in\X,
\quad\hbox{for every}\quad
f\in\X^*\!.
$$
All these equivalent formulations will be denoted by
$$
T_n\!\wconv T.
$$
\vskip-5pt\noi
Therefore, by Proposition~2.4,
$$
\hbox{\kern9.5pt\it $\{T_n\}\kern-1pt$ converges weakly
$\;\;$if and only if$\;\;$
$\,T_n\!\wconv T$ for some\/ $T\kern-1pt\in\LX$},                \eqno{(2.1)}
$$
and we will use the above equivalence freely throughout the text.

%%%%%%%%%%%%%%%%%%%%%%%%%%%  PROPOSITION 2.5
\vskip6pt\noi
{\bf Proposition 2.5.}\
{\it If\/ $\X$ is a Banach space and\/ $\{T_n\}$ is a\/ $\BX$-valued
sequence,~then\/ ${T_n\!\wconv T}$ implies\/ ${T\kern-1pt\in\BX}$ and\/
${\sup_n\|T_n\|\kern-1pt<\kern-1pt\infty}.$ $($So if\/ $\kern-1ptT_n\kern-1pt$
are bounded and act on a Banach space\/ $\X$, then the linear transformation\/
$T\kern-1pt$ in Proposition 2.4(b)~is~bounded}\/).

\proof
By Definition 2.3, Remark 2.2 ensures that
${\sup_n\|T_nx\|\kern-1pt<\kern-1pt\infty}$ for \hbox{every} ${x\in\X}.$
Thus if $\{T_n\}$ is a $\BX$-valued sequence and $\X$ is a Banach space, then
the \hbox{Banach}--Steinhaus Theorem ensures that
${\sup_n\|T_n\|\kern-1pt<\kern-1pt\infty}.$ Since $f(Tx)=\lim_nf(T_nx)$, we
get
$|f(Tx)|={\lim_n|f(T_nx)|}\le{\sup_n|f(T_nx)|}\le{\|f\|\sup_n\|T_n\|\,\|x\|}$
for every ${x\in\X}$, for every ${f\in\X^*}\!.$ Therefore, as a consequence
of the Hahn--\hbox{Banach} Theorem,
$\|Tx\|={\sup_{f\in\X^{^*}\!\!,\,\|f\|=1}|f(Tx)|}\le{\sup_n\|T_n\|\,\|x\|}$,
so that $T$ is bounded.                                                 \qed

\vskip9pt
If $T_n$ and $T$ lie in $\BX$, then ${f\kern1ptT_n}$ and ${f\kern1pt T}$ lie
in $\X^*$ for \hbox{every ${f\in\X^*}\!$,} so that the assertion
${f(T_nx)\to f(Tx)}$ for every ${x\in\X}$, for every ${f\kern-1pt\in\X^*}$
(i.e., ${T_n\!\wconv T}$), means that each $\X^*$-valued sequence
$\{f\kern1pt T_n\}$ converges pointwise to ${f\kern1pt T}$; that is,
$$
(f\kern1pt T_n)(x)\to(f\kern1pt T)(x)
\quad\hbox{for every}\quad
x\in\X,
\quad\hbox{for every}\quad
f\in\X^*\!.
$$
\vskip-1pt

\vskip5pt\noi
Next take any ${T\kern-1pt\in\LX}$ on a normed space $\X$ and consider its
{\it power sequence}\/ $\{T^n\}.$ Suppose the $\LX$-valued sequence $\{T^n\}$
converges weakly; that is, suppose
$$
T^n\!\wconv P
$$
for some ${P\in\LX}$ (according to Proposition 2.4)$.$ The special case where
$P$ is the null operator $O$ plays a special role$.$ We say that
${T\kern-1pt\in\LX}$ is {\it weakly stable}\/ if the $\LX$-valued power
sequence $\kern-.5pt\{T^n\}\kern-.5pt$ converges weakly to the null
\hbox{operator, denoted by}
$$
T^n\!\wconv O.
$$
\vskip-1pt

%%%%%%%%%%%%%%%%%%%%%%%%%%%  PROPOSITION 2.6
\vskip5pt\noi
{\bf Proposition 2.6.}\
{\it If\/ ${T\in\LX}$ on a normed space\/ $\X$, then
the following assertions are equivalent}\/.
\begin{description}
\item{$\kern-7pt$\rm(i)$\kern-1pt$}
{\it The power sequence\/ $\{T^n\}$ converges weakly}\/,
\vskip4pt
\item{$\kern-9pt$\rm(ii)}
${T^n\!\wconv P}$ {\it for some}\/ ${P\in\LX}$.
\end{description}
{\it In this case, and if\/ $T$ is bounded $($i.e., if\/
{\rm(ii)} holds and\/ ${T\kern-1pt\in\BX}\kern.5pt)$, then}
\begin{description}
\item{$\kern-9pt$\rm(a)}
$P$ {\it is a projection with\/ ${\R(P)=\N(I-T)},\;$ and}
\vskip4pt
\item{$\kern-9pt$\rm(b)}
${T^nP=PT^k}$ {\it for every pair of nonnegative integers\/ $n$ and}\/ $k$.
\vskip1pt\noi
\item{$\kern3pt$}
{\it In other words,\/ $\;{TP=T^nP=P=P^2=P\kern1ptT^k=P\kern1ptT}$ for
every}\/ ${n,k\ge2}$.
\end{description}
{\it Furthermore, still under the assumption that}\/ ${T\kern-1pt\in\BX}$ on
a normed space\/ $\X$,
\begin{description}
\item{$\kern-8pt$\rm(c)}
{\it ${T^n\!\wconv O}$ if and only if\/ ${T^n\!\wconv P}$ and}\/
${\N(I-T)=\0}$.
\end{description}
\begin{description}
\item{$\kern3pt$}
{\it Equivalently}\/,
$\;{T\kern-1pt\in\BX}$ {\it is weakly stable if and only if\/ $\{T^n\}$
converges weakly and\/ 1 is not an eigenvalue of}\/ $T$.
\end{description}

\proof
Take any ${T\kern-1pt\in\LX}$ and consider the power sequence $\{T^n\}.$ That
(i) and (ii) are equivalent follows from Proposition 2.4$.$ Suppose (i) or
(ii) holds so that ${T^n\!\wconv P}$ for some ${P\in\LX}.$

\vskip5pt\noi
First we prove (b); that is, if ${T\kern-1pt\in\BX}$, then
$$
T^nP=PT^k
\quad\hbox{for every}\quad
k\ge0
\;\;\hbox{and}\;\;
n\ge0.
$$
Indeed, if ${T\kern-1pt\in\LX}$ and ${f(T^nx)\to f(Px)}$ for every ${x\in\X}$,
for every ${f\in\X^*}\!$, then
\vskip5pt\noi
$$
{f(Px)\kern-1pt=\kern-1pt
{\lim}_n(T^{n+k}x)\kern-1pt=\kern-1pt
{\lim}_nf(T^nT^kx)\kern-1pt=\kern-1pt
{\lim}_nf\big(T^n(T^kx)\big)\kern-1pt=\kern-1ptf
(PT^kx)}
$$
\vskip3pt\noi
for every ${x\in\X}$, every ${f\kern-1pt\in\X^*}\!$, and every integer
${k\ge0}.$ If, in addition, ${T\kern-1pt\in\BX}$ (i.e., if the linear $T$ is
bounded), then ${f\kern.5ptT^n}$ lies in $\X^*$ for every
${f\kern-1pt\in\X^*}\kern-1pt$ and every integer ${n\ge0}$, so that
${\lim_k(f\kern1ptT^n)(T^kx)=(f\kern1ptT^n)(Px)}$ (because
${T^k\!\wconv P}).$ Hence
$$
{f(Px)\kern-1pt=\kern-1pt
{\lim}_k(T^{n+k}\kern-.5ptx)\kern-1pt=\kern-1pt
{\lim}_kf(T^nT^k\kern-.5ptx)\kern-1pt=\kern-1pt
{\lim}_k(f\kern.5ptT^n)(T^k\kern-.5ptx)\kern-1pt=\kern-1pt
(f\kern.5ptT^n)(Px)\kern-1pt=\kern-2pt
f(T^nPx)}
$$
\vskip3pt\noi
for every ${x\in\X}$, every ${f\kern-1pt\in\X^*}\!$, and every integer
${n\ge0}.$ Therefore, for each ${x\in\X}$,
$$
f(PT^kx)=f(T^nPx)
\quad\hbox{for every}\quad
f\kern-1pt\in\X^*.
$$
Then ${f(PT^kx-T^nPx)=0}$ for all ${f\in\X^*\!}$, so that ${PT^kx=T^nPx}$ for
every ${x\in\X}.$ Thus we get (b)$.$

\vskip5pt\noi
Proof of (a)$.$ With ${k=0}$ in (b), it follows that $P$ is a projection$:$
$$
Px=T^nPx\wconv PPx=P^2x
\quad\hbox{for every}\quad
x\in\X.
$$
To complete the proof of (a), take any ${x\in\R(P)}$ so that ${Px=x}$ (because
$P$ is a pro\-jection)$.$ Hence ${Tx=TPx}.$ But ${P=TP}$ by (b) (with ${n=0}$
and ${k=1}$) so that $Tx={Px=x}.$ Thus ${x\in\N(I-T)}.$ So
${\R(P)\sse\N(I-T)}.$ Conversely, if ${x\in\N(I-T)}$ (i.e., if ${Tx=x}$),
then ${T^nx=x}$ for every ${n\ge0}$, and hence ${Px=x}$ (as ${T^n\wconv P}$)
so that ${x\in\R(P)}.$ Thus ${\N(I-T)\sse\R(P)}.$ Therefore,
$$
\R(P)=\N(I-T).
$$
Proof of (c). Immediate: ${P=O\!\iff\!\R(P)=\0\!\iff\!\N(I-T)=\0.}$      \qed

\vskip5pt
Part of Proposition 2.6 has been shown in \cite[Theorem 3.2, Remark 3.3]{JJKS}
for operators acting on a Hilbert space, and in \cite[Theorem 1]{KV1} on a
separable Hilbert space$.$ The result below follows at once from
Proposition 2.5 and Proposition 2.6.

%%%%%%%%%%%%%%%%%%%%%%%%%%%  COROLLARY 2.7
\vskip5pt\noi
{\bf Corollary 2.7.}\
{\it If\/ ${T^n\!\wconv P}$ for\/ ${T\kern-1pt\in\BX}$ on a Banach space\/
$\X$, then\/ $P$ is a con\-tinuous projection $($i.e., ${P\in\BX}\kern.5pt)$
onto}\/ $\N({I-T})$.

\vskip5pt
If $\H$ is a Hilbert space with inner product ${\<\,\cdot\kern2pt;\cdot\,\>}$,
then by the Riesz Representation Theorem for Hilbert spaces, the property of
Definition 2.1, viz., there is a ${z\in\H}$ for which
${f(x_n)\to f(z)}$ for every ${f\kern-1pt\in\H^*}\!$, is equivalent to saying
that
$$
\hbox{there exists ${z\in\H}$ for which ${\<x_n\,;y\>\to\<z\,;y\>}$ for every
${y\in\H}.$}
$$
So weak convergence for an $\H$-valued sequence $\{x_n\}$ as in Definition 2.1
is rewritten as above$.$ Thus, according to Definition 2.3, an
$\Le[\H]$-valued sequence $\{T_n\}$ {\it converges weakly}\/ if for each
vector ${x\in\H}$ the $\H$-valued sequence $\{T_nx\}$ converges weakly, which
means that the scalar-valued sequence
$$
\hbox{
$\{\<(T_nx\,;y\>\}$ converges in $(\FF,|\cdot|)$ for every $x\in\H$, for
every $y\in\H$}.
$$
Now let $\kern-1pt\{T_n\}\kern-1pt$ be a $\BH$-valued sequence$.$ By
Proposition 2.5, this holds if and~only~if
\vskip4pt\noi
$$
\hbox{for each $x\in\H$, $\;\<T_nx\,;y\>\to\<Tx\,;y\>$ for every $y\in\H$}
$$
\vskip1pt\noi
(i.e., ${T_n\wconv T}$) for some operator ${T\kern-1pt\in\BH}$, which means
\vskip4pt\noi
$$
\hbox{$|\<(T_n-T)x\,;y\>|\to0$ for every $x\in\H$, for every $y\in\H$};
$$
\vskip1pt\noi
equivalently, if $\H$ is {\it complex}\/ (or $\kern-1pt\H$ real and
$T\kern-1pt$ self-adjoint; cf.\ polarisation identity),
\vskip4pt\noi
$$
\hbox{$|\<(T_n-T)x\,;x\>|\to0$ for every $x\in\H$}.
$$
\vskip1pt\noi
In such a Hilbert-space setting, weak convergence ${x_n\!\wconv x}$ of an
$\H$-valued sequence as in Definition 2.1 coincides with convergence in the
weak topology on $\H$, and the above definition of weak convergence
${T_n\!\wconv T}$ of a $\BH$-valued sequence as in Definition 2.3 coincides
with convergence in the weak operator topology on $\BH$ (see, e.g.,
\cite[Proposition IX.1.3]{Con} --- compare with the warning following
Definition~2.1).

%%%%%%%%%%%%%%%%%%%%%%%%%%%%%%%%%%%%%%%%%%%%%%%%%%%%%%%%%  SECTION 3
\section{Basic Notation and Terminology}

Throughout the paper, $\NN$, $\ZZ$, $\RR^{\kern-.5pt+}\kern-2.5pt,$ and
$\QQ$ will stand for the sets of positive~\hbox{integers}, of all integers, of
all nonnegative real numbers, and for the field of rational numbers.\ Let
$\DD$ and $\TT$ stand for the open unit disc and the unit circle in the
complex plane~$\CC$, and let $\X$ and $\Y$ be normed spaces$.$ We use the same
symbol ${\|\cdot\|}$ for the norm on $\X$ and for the induced uniform norm on
$\B[\X,\Y].$ A linear manifold $\M$ of $\X$ is non\-trivial if
${\0\ne\M\ne\X}.$ The closure of a linear manifold $\M$ will be denoted by
$\M^-\!.$ A closed linear manifold of a normed space $\X$ will be referred to
as a subspace~of~$\X$.

\vskip5pt
An operator $T\kern-1pt$ on a normed space $\X$ (i.e.,
${T\kern-1.5pt\in\kern-1pt\BX}$) is a contraction if ${\|T\|\le1}$, or
equivalently, if ${\|Tx\|\le\|x\|}$ for every ${x\in\X}\kern-1pt.$ It is a
strict contraction if ${\|T\|\kern-1pt<\kern-.5pt1}.$ In between there is the
notion of proper contraction, defined as ${\|Tx\|\kern-1pt<\kern-.5pt\|x\|}$
for \hbox{every} ${x\in\X}\kern-1pt.$ An isometry is an operator $T\kern-1pt$
such that ${\|Tx\|=\|x\|}$ for \hbox{every} ${x\in\X}\kern-1pt.$ An operator
$T\kern-1pt$ is power bounded if ${\sup_n\|T^n\|\kern-1pt<\kern-.5pt\infty}.$
(If $\X\kern-1pt$ is a Banach space, then power boundedness is equivalent to
${\sup_n\|T^nx\|\kern-1pt<\kern-.5pt\infty}$ for every ${x\in\X}\kern-1pt.)$ A
normed-space operator $T\kern-1pt$ is normaloid if ${\|T^n\|=\|T\|^n}\!$ for
all positive integers $n.$ A subspace $\M$ of $\X\kern-1pt$ is $T$-invariant
if ${T(\M)\sse\M}.$ If $\X\kern-1pt$ is a complex Banach space, then let
$\sigma(T)$ denote the spectrum of an operator $T$ on $\X\kern-1pt$, let
$\{\sigma_{\kern-1ptP}(T),\sigma_{\kern-1ptR}(T),\sigma_{\kern-.5ptC}(T)\}$
be the classical partition of the spectrum consisting of point spectrum
$\sigma_{\kern-1ptP}(T)$, residual spectrum $\sigma_{\kern-1ptR}(T)$, and
continuous spectrum $\sigma_{\kern-.5ptC}(T)$, and let
${r(T)=\sup_{\alpha\in\sigma(T)}|\alpha|}$ be the spectral radius of $T.$ The
\hbox{Gelfand}--\hbox{Beurling} formula says ${r(T)=\lim_n\|T^n\|^{1/n}}$
(where $\lim_n\|T^n\|^{1/n}$ exists even for normed-space operators)$.$ Recall
that ${r(T)^n=r(T^n)\le\|T^n\|\le\|T\|^n}$ for every positive integer $n$,
and note that $T\kern-1pt$ is normaloid if and only~if~${r(T)=\|T\|}$.

\vskip5pt
For any normed space $\X$, an $\X$-valued sequence ${\{x_n\}}$ converges
strongly (i.e., in the norm topology on $\X$) to ${x\in\X}$ if
${\|x_n\!-x\|\to0}.$ Notation$:$ ${x_n\!\to x}.$ It is clear that strong
convergence implies weak convergence (to the same and unique limit):
$$
x_n\!\to x
\;\quad\hbox{implies}\quad\;
x_n\!\wconv x.
$$
A $\BX$-valued sequence $\{T_n\}$ converges uniformly or strongly if (i) it
converges in the operator norm topology on $\BX$, or (ii) the $\X$-valued
sequence $\{T_nx\}$ converges in the norm topology on $\X$ for every
${x\in\X}$, respectively$.$ It is clear that if $\{T_n\}$ converges uniformly,
then it converges to an operator $T$ in $\BX.$ If $\X$ is a Banach space and
$\{T_n\}$ converges strongly, then it converges to an operator $T$ in $\BX.$
Thus uniform and strong convergences mean ${\|T_n\!-T\|\to0}$ and
${\|(T_n\!-T)x\|\to0}$ for every ${x\in\X}$ (notation$:$ ${T_n\!\uconv T}$ and
$\,{T_n\!\sconv T}$), respectively$.$ Again, it is clear that strong
convergence implies weak convergence, and these convergences are all to the
same limit (when they hold), so that
$$
T_n\!\uconv T
\quad\limply\quad
T_n\!\sconv T
\quad\limply\quad
T_n\!\wconv T.
$$
\vskip-1pt

\vskip5pt
Let $\H$ be a Hilbert space$.$ If $\M$ is a linear manifold of $\H$, then
$\M^\perp$ denotes~the orthogonal complement of $\M.$ Let $T^*\!$ in $\BH$ be
the (Hilbert-space) adjoint of $T\kern-1pt$ in $\BH.$ A subspace $\M$ of a
Hilbert space is reducing for an operator $T$ (or $\M$ reduces $T$) if $\M$
and $\M^\perp\!$ are both $T$-invariant or, equivalently, if $\M$ is invariant
for both $T$ and $T^*\!.$ A part of a Hilbert-space operator is a restriction
of it to a reducing sub\-space$.$ (Sometimes the term `part' of an operator
is defined as a restriction of the operator to an invariant subspace, but here
we define it a restriction to a reducing subspace$.$) The orthogonal direct
sum of spaces (or subspaces), and also of opera\-tors, will be denoted by
$\oplus.$ In such a Hilbert-space setting, an operator $T\kern-1pt$ is an
isometry if and only if ${T^*T\kern-1pt=\kern-1ptI}$, where $I$ stands for the
identity operator, and $T\kern-1pt$ is a coisometry if $T^*\!$ is an
isometry$.$ An operator $T$ is self-adjoint if ${T^*\!=T}$ and normal if
${TT^*\!=T^*T}$ (i.e., $T$ and $T^*\!$ commute)$.$ A unitary operator is an
invertible isome\-try, that is, a surjective isometry; equivalently, a normal
isometry, which means an isometry and a coisometry (i.e.,
${TT^*\!=T^*T\!=I}).$ An operator $T$ is quasinormal if
${(TT^*\!-T^*T)\kern1pt T\!=O}$ (i.e., if $T$ commutes with ${T^*T}$),
subnormal if it is the restriction of a normal operator to an invariant
subspace (i.e., if it has a normal extension), and hyponormal if
${TT^*\!\le T^*T}.\kern-1pt$ These classes are related by proper
inclusions~as~follows:
\vskip-2pt\noi
$$
\vbox{\noi
{\narrower\narrower
{\it Unitary and Self-Adjoint are Normal; Normal is Quasinormal, which \\
is Subnormal, which in turn is Hyponormal, all being Normaloid}\/.
\vskip0pt}
\vskip-20pt
}                                                                \eqno{(3.1)}
$$
\vskip11pt\noi
The existence of quasinormal operators that are not subnormal (in particular,
of normal operators that are not subnormal) shows that the restriction of a
quasinormal operator (or of a normal operator) to an invariant subspace is not
necessarily quasinormal (or not necessarily normal)$.$ However, the
direct sum of two operators is quasinormal (or normal) if and only if both
operators are quasinormal (or normal), and so parts of a quasinormal operator
(or parts of a normal operator) are trivially quasinormal (normal)$.$
Differently from these two cases,
\vskip-2pt\noi
$$
\vbox{\noi
{\narrower\narrower
{\it the restriction of a subnormal or of a hyponormal operator to \\
an invariant subspace is subnormal or hyponormal, respectively}\/.
\vskip0pt}
\vskip-20pt
}                                                                \eqno{(3.2)}
$$
\vskip11pt\noi
Indeed, $T$ is hyponormal if and only if $\|T^*x\|\le\|Tx\|$ for every
${x\in\H}.$ Thus, if $\M$ is $T$-invariant,
so that $(T|_\M)^*\!=E\,T^*|_\M$, where $E$ is the orthogonal projection onto
$\M$, then $\|(T|_\M)^*u\|\le\|T^*u\|\le\|Tu\|=\|T|_\M u\|$ for every
${u\in\!\M}.$ For the subnormal case, if $\M$ is $T$-invariant, then
${T\!=\!\big(\smallmatrix{T' & X' \cr
                          O' & Y' \cr}\big)\!:\!
                          \H=\!\M\oplus\M^\perp\!\to\H=\!\M\oplus\M^\perp\!}$
with ${T'\!=T|_\M\!:\!\M\to\!\M}.$ If $T$ is subnormal, then
${N\!=\!\big(\smallmatrix{T & X \cr
                          O & Y \cr}\big)\!:\!\H\oplus\K\to\H\oplus\K}$
for a normal operator $N.$ Hence
${N\!=\!\big(\smallmatrix{T'  & X'' \cr
                          O'' & Y'' \cr}\big)\!:\!
             \M\oplus(\M^\perp\!\oplus\K)\to\!\M\oplus(\M^\perp\!\oplus\K)}$.)
So
\vskip-2pt\noi
$$
\vbox{\noi
{\narrower\narrower
{\it parts of normal, quasinormal, subnormal, or hyponormal operators \\
are normal, quasinormal, subnormal, or hyponormal, respectively}\/.
\vskip0pt}
\vskip-20pt
}                                                                \eqno{(3.3)}
$$
\vskip11pt\noi
It is clear that every isometry is quasinormal$.$ A unilateral shift (which
is a nonnormal isometry) and a bilateral shift (which is a normal isometry,
i.e., unitary), of any multiplicity, will be denoted by $S_+$ and $S$,
respectively$.$

\vskip5pt
Let ${\<\,\cdot\kern2pt;\cdot\,\>}$ denote the inner product on a Hilbert
space $\H$, and let ${\|\cdot\|}$ be the norm on $\H$ induced by the inner
product$.$ Recall that an operator $T$ on a Hilb ert space $\H$ is self-adjoint
if ${T^*\!=T}$ (equivalently, if ${\<Tx\,;x\>\in\RR}$ for every ${x\in\H}$
whenever $\H$ is complex), and nonnegative if it is self-adjoint and
${\<Tx\,;x\>\ge0}$ for every ${x\in\H}.$ Also recall from the previous section
that an $\H$-valued sequence ${\{x_n\}}$ converges weakly to ${z\in\H}$ (i.e.,
${x_n\wconv z}$) if (and only if) ${\<x_n\,;y\>\to\<z\,;y\>}$ for every
${y\in\H}).$ In particular, in such a Hilbert-space setting, ${x_n\!\wconv z}$
implies ${\<x_n\,;z\>\to\<z\,;z\>=\|z\|^2}.$ Since
${\|x_n\!-z\|^2}=\|x_n\|^2-{2\,\hbox{Real}\,\<x_n\,;z\>+\|z\|^2}$, we get
$$
x_n\to z
\;\quad\hbox{if $\,$and $\,$only $\,$if}\quad\;
x_n\wconv z
\,\;\;\hbox{and}\;\;\,
\|x_n\|\to\|z\|.
$$
A Radon--Riesz space is a normed space $\X$ for which an $\X$-valued sequence
$\{x_n\}$~converges strongly if and only if it converges weakly (in the sense
of Definition 2.1)~and the sequence of norms $\{\|x_n\|\}$ converges to the
norm of the limit \cite[\hbox{Definition} 2.5.26]{Meg}$.$ The above-displayed
equivalence is called the Radon--Riesz property, and~we have just seen that
every Hilbert space is a Radon--Riesz space
\cite[{\hbox{Problem}~20,~p.13}]{Hal2}$.\!$

\vskip5pt
It is worth noticing that, for Hilbert-space operators,
$$
{\kern-55pt}
\hbox{\it weak convergence is preserved under the adjoint operation\/,}
                                                                \eqno{(3.4)}
$$
that is, $\!\{T_n\}$ converges weakly to $T$ if and only if
$\{{T_n}^{\kern-4pt*}\}$ converges weakly to $T^*\!$, as is trivially
verified once ${|\<(T_n\!-T)x\,;y\>|}={|\<({T_n}^{\kern-4pt*}\!-T^*)y\,;x\>|}$
for every ${x,y\in\H}.\kern-1pt$ {\it The same holds for uniform
convergence}\/ since $\|T_n-T\|=\|{T_n}^{\kern-4pt*}-T^*\|.$ However, the
power sequence of a unilateral shift shows that {\it this fails for strong
convergence}\/.

%%%%%%%%%%%%%%%%%%%%%%%%%%%%%%%%%%%%%%%%%%%%%%%%%%%%%%%%%  SECTION 4
\section{Notation and Terminology of Operator Stability}

Stability of operators means convergence of power sequences to the null
operator$.$ The nomenclature is usual, and its origin alludes to
asymptotic stability of discrete-time linear dynamical systems$.$ Indeed,
let $T$ be an operator acting on a normed space $\X$, and consider a discrete,
time-invariant, free, linear dynamical system modelled by the following
autonomous homogeneous difference equation,
$$
x_{n+1}=Tx_n
\quad\hbox{for every nonnegative integer $n$},
$$
for an arbitrary initial condition ${x_0\kern-1pt=x}$ in $\X$, whose solution
is ${x_n=T^nx}.\kern-1pt$ The linear system ${x_{n+1}=Tx_n}$ is
{\it asymptotically stable}\/ if the solution $\{x_n\}$ converges to zero for
every initial condition $x_0$ in $\X.$ Equivalently, if the $\X$-valued
sequence $\{T^nx\}$ converges to zero for every $x$ in $\X.$ According to the
nature of such a convergence, we have different notions of stability$.\!$
Therefore, an operator $T$ in $\BX$ is {\it uniformly stable}\/ if
${\sup_{x\ne0}\|T^nx\|/\|x\|\to0}$ (i.e., ${\|T^n\|\to0}$), {\it strongly
stable}\/ if ${\|T^nx\|\to0}$ (i.e., ${T^nx\to0}$) for every $x$ in $\X$, and
{\it weakly stable}\/ if $f{(T^nx)\to0}$ for every $x$ in $\X$, for every $f$
in $\X^*\!$ (i.e., ${T^nx\wconv0}$ for every $x$ in $\X$ as in Definition 2.1
--- in a Hilbert-space setting this means ${\<T^nx\,;y\>\to0}$ for every
${x,y}$ in $\H$ and, if $\H$ is complex, this means ${\<T^nx\,;x\>\to0}$ for
every $x$~in~$\H$).\ So an operator ${T\kern-1pt\in\BH}$ is {\it uniformly
stable}\/, {\it strongly stable}\/, or {\it weakly stable}\/ if\/
${T^n\!\uconv O}$, ${T^n\!\sconv O}$, or ${T^n\!\wconv O}$,~respectively.\
(Such a framework is standard --- see, e.g., \cite[Section 2]{KV1}; see also
\cite[Section 1]{JJKS}.)

\vskip5pt
Uniform stability is fully characterised$.$ A Banach-space operator is
uniformly~stable if and only if its spectral radius is less than one:
${\kern-1pt T^n\kern-2pt\uconv\kern-1pt O\kern-3pt
\iff\kern-3pt r(T)\kern-1.5pt<\kern-1pt1}$~(see,~e.g.,
\cite[Proposition 0.4]{MDOT}, \cite[Proposition II.1.3]{Eis}$\kern.5pt$).\
The spectral radius of a power-bound\-ed \hbox{Banach}-space operator is not
greater than one:\
${\sup_n\|T^n\|\kern-1pt<\kern-1pt\kern-1pt\infty
\limply r(T)\kern-1pt\le\kern-1pt1}.$
$\kern-1pt$This follows from the Gelfand--Beurling formula,
${r(T)=\lim_n\|T^n\|^{1/n}}\!$, and the converse fails:\ there exist a
diagonal operator $T\kern-1pt$ on $\ell_+^2$, the Hilbert space of square
summable scalar sequences, such that ${\|T^n\|=(\log n)^{1/2}}$ (see, e.g.,
\cite[Section III.4.A, p.66]{Bea}$\kern.5pt$), and hence
${r(T)^2=\lim_n(\log n)^{1/n}=1}.$ Since weak
stability of a \hbox{Banach}-space~operator implies power boundedness
(cf.\ Proposition~2.5), we get
\begin{eqnarray*}
r(T)<1\;\;
&\kern-6pt\iff\kern-6pt&
\;\;T^n\!\uconv O\;\;\limply                                               \\
T^n\!\sconv O\;\;\limply\;\;T^n\!\wconv O\;\;
&\kern-6pt\;\limply\kern-6pt&
\;\;{\sup}_n\|T^n\|<\infty\;\;\limply\;\;r(T)\le1,
\end{eqnarray*}
and the reverses of the above one-way implications fail$.$ However, on a
Hilbert space,
$$
{\kern-47pt}
\hbox{\it weak and strong stabilities coincide for self-adjoint operators}
                                                                \eqno{(4.1)}
$$
(i.e., ${T^{*n}T^n\!\wconv\kern-1ptO}$ means ${\|T^nx\|^2\!\to0}$ so that
${T^n\!\wconv\kern-1ptO}$ implies ${T^n\!\sconv O}$ if ${{T^*}\!=T\!}.)$
This fails for unitary operators, which are isometries, and so never strongly
stable.\ As power boundedness implies a spectral radius not greater than one,
we get:
$$
\hbox{\it if a power bounded operator is normaloid, then it is a contraction.}
$$
Thus the above string of implications leads to the following useful results.
$$
\vbox{\noi
{\narrower\narrower
{\it Suppose\/ an operator on a complex Banach space is\/ {\rm normaloid}.\
Then
\vskip5pt\noi
${\kern13pt}${\rm(i)}
it is uniformly stable if and only if it is a strict contraction,
\vskip1.5pt\noi
${\kern12pt}${\rm(ii)}
it is power bounded if and only if it is a contraction, $\;$and so
\vskip1.5pt\noi
${\kern10pt}${\rm(iii)}
if it is weakly stable, then it is a
contraction}\/.
\vskip0pt}
\vskip-35pt
}                                                               \eqno{(4.2)}
$$
\vskip26pt\noi
Moreover, (see, e.g., \cite[Proposition 2.O]{ST2}\kern.5pt)
\vskip0pt\noi
$$
\vbox{\noi
{\narrower\narrower
{\it the concepts of weak, strong, and uniform stabilities coincide for \\
a compact operator on a complex Banach space}\/.
\vskip0pt}
\vskip-20pt
}                                                               \eqno{(4.3)}
$$
\vskip6pt

\vskip5pt
Resuming the Hilbert-space setting of Section 3 in light of stability, the
power~sequence $\{T^n\}$ of a Hilbert-space operator ${T\kern-1pt\in\BH}$
converges weakly if there~is~an~operator ${P\kern-1pt\in\BH}$ for which
${T^nx\wconv Px}$ for every ${x\in\H}$, that is, ${\<T^nx\,;y\>\to\<Px\,;y\>}$
for every ${x,y\in\H}.$ As we saw at the end of Section 2, this in turn is
equivalent to saying that ${\<T^nx\,;x\>\to\<Px\,;x\>}$ for every ${x\in\H}$
if $\H$ is a complex Hilbert space$.$~No\-tation:\
$\kern-1pt{T^n\!\wconv\kern-1pt P}\kern-1pt.$ If $P$ is the null operator $O$,
then $T$ is weakly stable, which~means
$$
\<T^nx\,;y\>\to0
\hbox{ for every }
x,y\in\H;
\quad\hbox{equivalently,}\quad
\<T^nx\,;x\>\to0
\hbox{ for every }
x\in\H
$$
if the Hilbert space is complex$.$ The Radon--Riesz property says that
$$
T^n\!\sconv P
\quad\hbox{if and only if}\quad
T^n\!\wconv P
\,\;\hbox{and}\;\,
\|T^nx\|\to\|Px\|
\,\;\hbox{for every}\;\,
x\in\H.
$$
Such an equivalence shows that a prototype of weakly stable but not strongly
stable operators are weakly stable isometries, which are characterised below.

\vskip5pt
The celebrated {\small \sc Nagy--Foia\c s--Langer Decomposition} for
contractions \cite{NF1}, \cite{Lan} (see also \cite[Theorem 3.2]{NF} or
\cite[Theorem 5.1]{MDOT}$\kern.5pt$) says that {\it every Hilbert-space
contraction\/ $T$ can be uniquely decomposed as an orthogonal direct sum of
two parts\/, a completely nonunitary contraction\/ $C$ and a
unitary operator\/ $U\kern-1pt$}\/$,$
$$
T\kern-1pt=C\oplus U.
$$
(A Hilbert-space contraction is completely nonunitary if the restriction of it
to every reducing subspace is not unitary.) The {\small \sc von Neumann--Wold
Decomposition} (see, e.g., \cite[Theorem 1.1]{NF} or
\cite[Corollary 5.6]{MDOT}$\kern.5pt$) can be viewed as a special case of the
Nagy--Foia\c s--Langer decomposition, and it says that {\it every
Hilbert-space isometry\/ $V\!$ has a unique decomposition as an orthogonal
direct sum of a unilateral shift\/ $($of some multiplicity\/$)$ $S_+\kern-1pt$
and a unitary operator\/ $U\kern-1pt$},
$$
V\kern-1pt=S_+\oplus U.
$$
Observe that any of the parts in the above two decompositions may be
missing$.$ As the von Neumann--Wold decomposition suggests (see also, e.g.,
\cite[\hbox{Lemma 5.4}]{MDOT}$\kern.5pt$),
$$
{\kern-35pt}
\hbox{\it a unilateral shift is precisely a completely nonunitary isometry}\/.
                                                                \eqno{(4.4)}
$$
Also, recall that
$$
\vbox{\noi
{\narrower\narrower
{\it
--- shifts $(\kern-1pt$unilateral or bilateral, $\kern-1pt$of any
multiplicity$)$ are weakly stable, \\
--- unilateral shifts are weakly stable isometries, \\
--- bilateral shifts are weakly stable unitaries}\/.
\vskip0pt}
\vskip-26pt
}                                                              \eqno{(4.5)}
$$
\vskip14pt

\vskip5pt
According to (3.4), weak stability is preserved under the adjoint operation
for~all Hilbert-space operators.\ This fails for strong stability.\ Indeed, a
unilateral shift is a (nonnormal) isometry with a strongly stable adjoint.\
However,
\vskip-2pt\noi
$$
\vbox{\noi
{\narrower\narrower
{\it strong stability is preserved under the adjoint operation for \\
normal operators}\/.
\vskip0pt}
\vskip-20pt
}                                                               \eqno{(4.6)}
$$
\vskip10pt\noi
(Reason$:$ ${T\!\in\kern-1pt\BH}$ is normal
$\!\iff\!\|T^{*n}x\|\kern-1pt=\kern-1pt\|T^nx\|$ for every
${x\kern-1pt\in\kern-1pt\H}$ and \hbox{every}~${n\kern-1pt\ge\kern-1pt0}.$)

%%%%%%%%%%%%%%%%%%%%%%%%%%%  REMARK 4.1
\vskip5pt\noi
{\bf Remark 4.1.}\
A suitable starting point for a review on weak stability is the elemen\-tary
result saying that weak stability travels well from parts to direct sums and
back.\ Let ${\K\oplus\Le}$ be the orthogonal direct sum of Hilbert spaces
$\K$ and $\Le$ equipped with its natural inner product,
${\<(x,y)\,;(w,z)\>_{\K\oplus\Le}}
\kern-1pt=\kern-1pt{\<x\,;w\>_\K}+{\<y\,;z\>_\Le}$
for every ${(x,y),(w,z)\kern-1pt\in\kern-1pt\K\oplus\Le}.$ Let
${A\oplus B=\big(\smallmatrix{A & O \cr
                              O & B \cr}\big)}$
and
${P\oplus E=\big(\smallmatrix{P & O \cr
                              O & E \cr}\big)}$
on ${\K\oplus\Le}$ be the direct sums of operators ${A,P\in\BK}$ and
${B,E\in\BL}.$ It is readily verified that
${(A\oplus B)^n}\kern-1pt={A^n\oplus B^n}\kern-1pt$ for every integer
${n\ge0}$, and
$$
(A\oplus B)^n\kern-1pt\wconv P\oplus E
\quad\hbox{if and only if}\quad
A^n\kern-1pt\wconv P
\,\;\hbox{and}\;\,
B^n\kern-1pt\wconv E.
$$
(This is naturally extended to multiple direct sums such as
${\bigoplus_{i=1}^mA_i}.)$ In particular,
$$
(A\oplus B)^n\kern-1pt\wconv O
\quad\hbox{if and only if}\quad
A^n\kern-1pt\wconv O
\;\;\hbox{and}\;\;
B^n\kern-1pt\wconv O.
$$
It is clear that the same argument holds for weak, strong, and uniform
convergences in general (not necessarily convergences of power sequences).

\vskip5pt
Therefore, since every unilateral shift is weakly stable, it follows by the
von \hbox{Neumann}--Wold decomposition that (if the unitary part is present)
\vskip-2pt\noi
$$
\vbox{\noi
{\narrower\narrower
{\it an isometry is weakly stable if and only if its unitary part is \\
weakly stable}\/.
\vskip0pt}
\vskip-20pt
}                                                               \eqno{(4.7)}
$$
\vskip14pt

%%%%%%%%%%%%%%%%%%%%%%%%%%%%%%%%%%%%%%%%%%%%%%%%%%%%%%%%%  SECTION 5
\section{Weakly Stable Operators}

Weakly stable operators have been discussed in \cite[Section 8.2]{MDOT} and
\cite[Section II.3]{Eis}, and also on a few pages of some books dealing
with functional \hbox{analysis}/oper\-ator theory in general$.$ We begin our
general exposition with the {\small \sc Foguel Decomposition}
\cite[Theorem 1.1]{Fog1} for Hilbert-space contractions, which decompose them
into a direct sum of a weakly stable contraction and a unitary operator.

%%%%%%%%%%%%%%%%%%%%%%%%%%%  PROPOSITION 5.1
\vskip5pt\noi
{\bf Proposition 5.1} \cite{Fog1} (1963).\
{\it Let\/ $T\kern-1pt$ be an arbitrary contraction on a Hilbert space\/
$\H.\!$ The set\/ $\Z(T)$ of all weakly stable vectors\/ $x$ for\/
$T\kern-1pt$,
$$
\Z(T)=\big\{x\in\H\!:T^nx\wconv0\big\},
$$
is a reducing subspace for\/ $T\kern-1pt$, and $T\kern-1pt$ is decomposed on\/
${\H=\Z(T)\oplus\Z(T)^\perp}$ as an orthogonal direct sum of a weakly stable
contraction and a unitary operator,
$$
T=Z\oplus U,
$$
 where\/ $Z=T|_{\Z(T)}$ is a weakly stable contraction and\/
$U=T|_{\Z(T)^\perp}$ is unitary}\/.

\vskip5pt
Unlike Nagy--Foia\c s--Langer decomposition, Foguel decomposition is not
unique; there are weakly stable unitary operators (e.g., bilateral shifts).\
So the part $Z$ of a Foguel decomposition of a contraction $T$ may not be
completely nonunitary.\ But Nagy--Foia\c s--Langer and Foguel decompositions
together ensure the following important result (see \cite[Corollary to
Theorem 6.III]{Fil}; see also \cite[\hbox{Corollary 7.4}]{MDOT}$\kern.5pt$).

%%%%%%%%%%%%%%%%%%%%%%%%%%%  COROLLARY 5.2
\vskip5pt\noi
{\bf Corollary 5.2} \cite{Fil} (1970).\
{\it $\kern-1pt$Every completely nonunitary contraction is weakly
\hbox{stable}}\/.

\vskip5pt
Thus Nagy--Foia\c s--Langer decomposition gives a contraction counterpart to
(4.7):
\vskip-3pt\noi
$$
\vbox{\noi
{\narrower\narrower
{\it a contraction is weakly stable if and only if its unitary part is \\
weakly stable}\/.
\vskip0pt}
\vskip-20pt
}                                                               \eqno{(5.1)}
$$
\vskip5pt

\vskip5pt
Strict contractions are trivially uniformly stable (since ${r(T)\le\|T\|}$),
and so they are clearly weakly stable$.$ Plain contractions are power bounded
but not necessari\-ly weakly stable (e.g., the identity)$.$ But proper
contractions are completely nonunitary (since ${\|Tx\|<\|x\|}$ for every
nonzero ${x\in\X}).$ Thus, Foguel decomposition, through Cor\-ollary 5.2,
ensures the next result from \cite[Proposition 2.2]{KL}.

%%%%%%%%%%%%%%%%%%%%%%%%%%%  COROLLARY 5.3
\vskip5pt\noi
{\bf Corollary 5.3} \cite{KL} (2001).\
{\it Every Hilbert-space proper contraction is weakly stable}\/.

\vskip5pt
It is clear that the converses of Corollaries 5.2 and 5.3 fail.\

\vskip5pt
If $T\kern-1pt$ is a hyponormal contraction (i.e., if $T^*\!$ is a
cohyponormal contraction),~then $T^*\!$ is uniquely decomposed as an
orthogonal direct sum of two parts:\ a strongly stable contraction and a
unitary operator.\ Indeed, if $T$ is a hyponormal contraction, then
\cite[Theorem 1 and Example 2]{KVP},
$$
T^*\!=G\oplus U,
$$
where $G$ is a strongly stable cohyponormal contraction and $U$ is a unitary
operator.\ Since the strongly stable $G$ is completely nonunitary and since
the unitary $\kern1ptU\kern-1pt$ has no strongly stable part, the
decomposition is unique, and any part of it may be~missing.\ In particular, if
$T$ is a normal contraction, then \cite[Corollary 1 and Example 4]{KVP}
$$
T=B\oplus U,
$$
where $B$ is a strongly stable normal contraction (with a strongly stable
adjoint --- cf.\ (4.6)$\kern.5pt$), and $U$ is unitary (see also
\cite[Section 5.3]{MDOT} and \cite[Section 2]{Kub2}$\kern.5pt$).\

\vskip5pt
On the other hand, according to (4.2), {\it if a Hilbert-space operator is
weakly stable and normaloid, then it is a contraction}\/.\ In particular, by
(3.1), {\it a weakly stable hyponormal\/ $($or normal\/$)$ operator is a
contraction}\/.\ So, recalling that weak convergence is preserved under the
adjoint operation as in (3.4), we get a characterisation of weakly stable
hyponormal and normal operators.\

%%%%%%%%%%%%%%%%%%%%%%%%%%%  PROPOSITION 5.4
\vskip5pt\noi
{\bf Proposition 5.4.}\
{\it Let\/ $\H$ be a Hilbert space}.\
\begin{description}
\item{$\kern-9pt$\rm(a)$\kern1pt$}
{\it A hyponormal operator\/ $T\kern-1pt$ on\/ $\H$ is weakly stable if and
only if\/ ${T^*\!=G\oplus U}\!$, where\/ $G$ is a strongly stable cohyponormal
contraction and\/ $U\!$ is a weakly stable unitary operator}\/.\
\vskip4pt
\item{$\kern-9pt$\rm(b)$\kern-0.5pt$}
{\it A normal operator\/ $T$ on\/ $\H$ is weakly stable if and only if\/
${T\kern-1pt=B\oplus U}\!$, where\/ $B$ is a strongly stable normal
contraction\/ $($with a strongly stable adjoint\/$)$ and\/ $U\!$ is a weakly
stable unitary operator}\/.\
\end{description}
\vskip-4pt\noi

\proof
(a)
Suppose $T$ is a weakly stable hyponormal operator$.$ As we saw above, the
results in (3.1) and (4.2) ensure that $T$ is a contraction$.$ Thus apply the
decomposition ${T^*\!=G\oplus U}$ for the cohyponormal contraction $T^*\!$,
where $G$ is a strongly stable cohyponormal contraction (because, according
to (3.3), parts of cohyponormal operators are again cohyponormal), and the
other part, the unitary operator $U\!$, is naturally cohyponormal.\ Remark 4.1
ensures that the unitary $U$ is weakly stable.\ The converse is trivial by
Remark 4.1 (as orthogonal direct sums of cohyponormal operators are
cohyponormal).\
\vskip5pt\noi
(b)
This is a particular case of (a) via the decomposition
${T\kern-1pt=\kern-1ptB\oplus U}\!$ for normal contractions (since normality
also travels well between parts and direct sum).                        \qed

\vskip5pt
The case for normal operators in Proposition 5.4(b) was proved in
\cite[Corollary 6.6]{JJKS} by using a different approach$.$ The decomposition
${T\kern-1pt=\kern-1ptB\oplus U}\!$ for normal con\-tractions and
Proposition 5.4(b) ensure that (compare with (4.1) and Corollary~5.2\/)
$$
\hbox{\kern-3pt
\it every completely nonunitary normal contraction is strongly
stable}\/.                                                      \eqno{(5.2)}
$$
\vskip-5pt

\vskip5pt
Recall from Proposition 2.6 and Corollary 2.7 that a power sequence of an
operator $T$ on a Hilbert space $\H$ (in general, on a Banach space $\X$)
converges weakly if and only if it converges weakly to an operator in $\BH.$
In other words,
$$
\hbox{\kern14pt\it $\{T^n\}\kern-1pt$ converges weakly
$\;\;$if and only if$\;\;$
$\,T^n\!\wconv P$ for some\/ $P\in\BH$},                        \eqno{(5.3)}
$$
and, again, we use the above equivalence freely throughout the text.

\vskip5pt
Proposition 5.5 below was given in \cite[Theorem 1]{KV1}$.$ The equivalence
(a)$\Leftrightarrow$(b) and item (e) were established there in a separable
Hilbert space, which have been extended here to normed and Banach spaces
in Propositions 2.6 and 2.7.

%%%%%%%%%%%%%%%%%%%%%%%%%%%  PROPOSITION 5.5
\vskip5pt\noi
{\bf Proposition 5.5} \cite{KV1} (1989).\
{\it If\/ $T\kern-1pt$ and\/ $\kern-.5ptP\kern-.5pt$ are operators on a
separable Hilbert space\/ $\H$ with an orthonormal basis\/
$\kern-.5pt\{e_k\}\kern-.5pt$, then the following assertions are equivalent\/.}
\begin{description}
\item{$\kern-6pt$\rm(a)}
${T^n\!\wconv P}$.
\vskip4pt
\item{$\kern-6pt$\rm(b)}
${PT=TP=P^2=P}$ {\it and}\/ ${(T^n\!-P)^n\!\wconv O}$.
\vskip4pt
\item{$\kern-6pt$\rm(c)}
$T$ {\it is power bounded and\/ ${\<T^ne_k;e_\ell\>\to\<Pe_k;e_\ell\>}$
for every integers}\/ $k,\ell$.
\end{description}
{\it Moreover, if any of the above equivalent conditions holds, then}
\begin{description}
\item{$\kern-6pt$\rm(d)}
${r(T)\le1}$,
\quad
${\sigma_{\kern-1ptR}(T)\sse\DD}$,
\quad
${\sigma_{\kern-1ptP}(T)\sse\DD}\cup\{1\}$,
\quad
{\it and}
\vskip4pt
\item{$\kern-6pt$\rm(e)}
$P=O\iff\sigma_{\kern-1ptP}(T)\sse\DD\iff1\not\in\sigma_{\kern-1ptP}(T)$.
\end{description}
\vskip-1pt

\vskip5pt
In fact, separability for $\H$ in Proposition 5.5 is necessary only for item
(c), which refers to a countable orthonormal basis $\{e_k\}.$ Actually, weak
stability always implies (i) point and residual spectra in the open unit disc,
(ii) spectrum in the closed unit disc, and (iii) spectral radius not greater
than one, as summarised next.

%%%%%%%%%%%%%%%%%%%%%%%%%%%  REMARK 5.6
\vskip5pt\noi
{\bf Remark 5.6.}\
Let\/ $T$ be an operator on a (complex) Hilbert space$.$ Then
\begin{description}
\item{$\kern-6pt$\rm(a)}
${T^n\!\wconv O\;\limply\;\sup_n\|T^n\|<\infty
\;\;\hbox{\it and}\;\;
\sigma_{\kern-1ptP}(T)\sse\DD}$,
\vskip4pt
\item{$\kern-6pt$\rm(b)}
${\sup_n\|T^n\|<\infty\;\limply\;\sigma_{\kern-1ptR}(T)\sse\DD
\;\;\hbox{\it and}\;\;
r(T)\le1}$,
\vskip4pt
\item{$\kern-6pt$\rm(c)}
${r(T)<1\iff T^n\!\uconv O\iff T^n\!\wconv O
\;\;\hbox{\it and}\;\;
\sigma_{\kern-.5ptC}(T)\cap\TT=\void}$.
\end{description}
\vskip-1pt
\vskip5pt\noi
Indeed, we had seen in Section 3 that part of (a), part of (b), and part of
(c),~namely,
${r(T)<1\,\Leftrightarrow\,T^n\!\uconv O\,\Rightarrow\,T^n\!\wconv O
\,\Rightarrow\,\sup_n\|T^n\|<\infty\,\Rightarrow\,r(T)\le1}$,
hold even in a \hbox{Banach} space$.$ The other part of (a) is trivial
(${T^n\!\wconv O\,\Rightarrow\,\sigma_{\kern-1pt}(T)\sse\DD}$), and the other
part of (b), viz.,
${\sup_n\|T^n\|<\infty\,\Rightarrow\,\sigma_{\kern-1pt R}(T)\sse\DD}$,
that is,
\vskip-3pt\noi
$$
\vbox{\noi
{\narrower\narrower
{\it the residual spectrum of a power-bounded operator lies in the         \\
open unit disc}\/,
\vskip0pt}
\vskip-20pt
}                                                                \eqno{(5.4)}
$$
\vskip10pt\noi
was shown in \cite[Theorem 4.1]{MK}.\ The remaining results, including the
equivalence,
\vskip-5pt\noi
$$
\vbox{\noi
{\narrower\narrower
{\it an operator is uniformly stable if and only if it is weakly stable    \\
and the continuous spectrum does not meet the unit circle}\/,
\vskip0pt}
\vskip-20pt
}                                                                \eqno{(5.5)}
$$
\vskip10pt\noi
can be found, for instance, in
\cite[Proposition 8.4 and Remark 8.6]{MDOT}$\kern.5pt$).
\vskip5pt
Notation$:$ we use a unique notation for the identity operator on any linear
space, so the same notation $I$ is used for the identity on $\H$ and also on
${\N(I\kern-1pt-T)}.$ Recall that ${\N(I\kern-1pt-T)=\0}$ means
${1\not\in\sigma_{\kern-1ptP}(T)}$ for any operator $T$.

\vskip5pt
Weak stability for unitary operators on an arbitrary Hilbert space was
recently characterised in \hbox{\cite[Lemma 3.6]{CJJS2}} as follows.

%%%%%%%%%%%%%%%%%%%%%%%%%%%  PROPOSITION 5.7
\vskip5pt\noi
{\bf Proposition 5.7} \cite{CJJS2} (2024).\
{\it Let\/ $U$ be a unitary operator on a Hilbert space}\/.
\begin{description}
\item{$\kern-6pt$\rm(a)}
{\it If\/ ${U^n\!\wconv P}$, then\/ $P$ is an orthogonal projection such
that}\/ ${\R(P)=\N(I\kern-1pt-U)}$.
\vskip4pt\noi
\item{$\kern-6pt$\rm(b)}
{\it $U$ is weakly stable if and only if\/ ${U^n\!\wconv P}$ and}\/
${\N(I\kern-1pt-U)=\0}$.
\end{description}
\vskip-1pt

\vskip5pt
The next result from \cite[Theorem 3.2, Remark 3.3]{JJKS} extends
Proposition 5.5 to an arbitrary Hilbert space and Proposition 5.7 to an
arbitrary Hilbert-space operator$.$ Moreover, it also gives a decomposition of
an arbitrary Hilbert-space operator whose power sequence is weakly convergent
into a direct sum of an identity operator and a weakly stable one$.$ Such a
decomposition comes about according to the following argument$.$ If $P$ is a
projection, then ${\N(P)=\R(I\kern-1pt-P)}$ and ${\R(P)=\N(I\kern-1pt-P)}$ are
algebraic complements of each other, leading to the \hbox{\it algebraic}\/
direct sum decomposition ${\H=\R(P)\oplus\R(I\kern-1pt-P)}.$ Also, if
${P\kern1pt T\!=TP}$, then $\R(P)$ and $\R({I\kern-1pt-P})$ are $T$-invariant,
and so, if ${\R(P)=\N(I\kern-1pt-T)}$, then we get the {\it algebraic}\/
direct sum decomposition $T\!={I\oplus L}$, with
${I\!=T|_{\R(P)}=T|_{\N(I-T)}}$ and ${L=T|_{\R(I-P)}}.$ (We~have used above
the same symbol $\oplus$ for {\it algebraic}\/ direct sums; not
necessarily~\hbox{orthogonal}).

%%%%%%%%%%%%%%%%%%%%%%%%%%%  PROPOSITION 5.8
\vskip5pt\noi
{\bf Proposition 5.8} \cite{JJKS} (2024).\
{\it Let\/ ${T\kern-1pt\in\BH}$ be an operator on a Hilbert space\/ $\H$.}
\begin{description}
\item{$\kern-9pt$\rm(a)$\kern-.5pt$}
{\it If\/ ${T^n\!\wconv P\in\BH}$, then $P$ is a projection onto\/ $\N({I-T})$
that commutes with\/ $T$ $($i.e.,\/ ${P=P^2}$ with\/ ${\R(P)=\N(I-T)}$ and
${P\kern1ptT\kern-1pt=TP}\/)$}\/.
\end{description}
{\it Moreover}\/,
\begin{description}
\item{$\kern-14pt$\rm(b$_1$)$\kern1pt$}
{\it ${T\kern-1pt=I\oplus L}$ on\/ ${\H=\R(P)\oplus\R(I-P)}$, where the
symbols $\oplus$ stand for algebraic\/ $($not necessarily orthogonal\/$)$
direct sums\/, $I$ stands for the identity operator on\/ $\R(P)$, and\/ $L$ is
a weakly stable operator on\/ $\R({I-P})$}\/.
\end{description}
{\it Conversely}\/,
\begin{description}
\item{$\kern-14pt$\rm(b$_2$)$\kern1.5pt$}
{\it if\/ ${P\in\BH}$ is a projection on\/ $\H$ and if\/ $T$ decomposes as
above\/ $($i.e., if\/~$T\kern-1pt=$ ${I\oplus L}$ on\/
${\H=\R(P)\oplus\R(I-P)}$ --- algebraic direct sums --- where\/ $I$ is the
identity and\/ $L$ is weakly stable\/$)$, then\/ ${T^n\!\wconv P}$}\/.
\end{description}
{\it Furthermore}\/,
\begin{description}
\item{$\kern-9pt$\rm(c)$\kern2pt$}
{\it $T\kern-1pt$ is weakly stable if and only if\/ $\{T^n\}$ converges
weakly and\/ ${\N(I-T)=\0}$}\/.
\end{description}
\vskip-1pt

\vskip5pt
Observe that the results from Proposition 2.4 to Corollary 2.7 show that
Propositions 5.5(a,b,e) and 5.8(a,c) hold for Banach-space operators.

\vskip5pt
From now on $\oplus$ will stand again for {\it orthogonal}\/ direct sum$.$
The proposition below is from \cite[Theorem 4.2 and Corollary 4.4]{JJKS}.

%%%%%%%%%%%%%%%%%%%%%%%%%%%  PROPOSITION 5.9
\vskip5pt\noi
{\bf Proposition 5.9} \cite{JJKS} (2024).\
{\it Let\/ ${T\kern-1pt\in\BH}$ be an operator on a Hilbert space\/ $\H$.
\vskip4pt\noi
{\rm (a)}
If\/ ${T^n\!\wconv P}$ $($so that ${P\in\BH}$ is a projection with\/
${\R(P)=\N(I-T)}$ according to Proposition 5.8$)$, then the following
assertions are equivalent}\/.
\begin{description}
\item{$\kern2pt$\rm(i)\kern3pt}
{\it The projection $P$ is orthogonal}\/.
\vskip3pt\noi
\item{$\kern1pt$\rm(ii)\kern3pt}
${\N(I-T)=\N(I-T^*)}$.
\vskip3pt\noi
\item{$\kern-.5pt$\rm(iii)\kern2pt}
${\N(I-T)}$ {\it reduces}\/ $T$.
\end{description}
\vskip3pt\noi
{\rm (b)}
{\it If\/ ${\N(I-T)}$ reduces\/ $T$, so that}
$$
T\kern-1pt=I\oplus L
\quad\hbox{on}\quad
\H=\N(I-T)\oplus\N(I-T)^\perp,
$$
{\it then\/ ${T^n\!\wconv P}$ if and only if\/ ${L^n\!\wconv O}.$ In this
case,\/ $P$ is the orthogonal projection\/ $P={I\oplus O}$ on\/
${\H=\N(I-T)\oplus\N(I-T)^\perp}$ $($so that}\/ ${\R(P)=\N(I-T)}\kern.5pt)$.

\vskip5pt
In fact, for an arbitrary Hilbert-space operator $\kern-1pt T$, if the
property ${\N(\alpha I-T)}\sse{\N(\overline\alpha I-T^*)}$ holds for some
${\alpha\in\FF}$, then ${\N(\alpha I-T)}$ reduces $T\kern-1pt$, and so
$T\kern-1pt$ may be (uniquely) decomposed as
$$
T\kern-1pt=\alpha I\oplus L
\quad\hbox{on}\quad
\H=\N(\alpha I-T)\oplus\N(\alpha I-T)^\perp\!,
$$
where ${\alpha I=T|_\N(I-T)}$, with $I$ denoting the identity operator on
${\N(\alpha I-T)}$, and ${L=T|_{\N(\alpha I-T)^\perp\!}}$ on
${\N(\alpha I-T)^\perp\!}.$ (It is clear that any part, namely, $\alpha I$ or
$L$, may be absent in the above direct sum decomposition).\ Recall that if
$T\kern-1pt$ is hyponormal, then ${\N(\alpha I-T)}\sse{\N(\alpha I-T^*)}$ for
every ${\alpha\in\FF}$ (see, e.g., \cite[Lemma 1.13]{ST2}$\kern.5pt$), and if
$T$ is a contraction, then ${\N(I-T)}={\N( I-T^*)}$ (see, e.g.,
\cite[Proposition I.3.1]{NF}$\kern.5pt).$ Therefore, in both cases,
$\N({I-T})$ reduces $T\kern-1pt$, and we get the following consequence of
Propositions 5.8 and 5.9, as in \cite[\hbox{Corollary 4.5}]{JJKS}.

%%%%%%%%%%%%%%%%%%%%%%%%%%%  COROLLARY 5.10
\vskip5pt\noi
{\bf Corollary 5.10} \cite{JJKS} (2024).\
{\it Let\/ ${T\kern-1pt\in\kern-.5pt\BH}$ be an operator on a Hilbert space}\/
$\H.$
\begin{description}
\item{$\kern-7pt$\rm(a)\kern-1pt}
{\it If\/ $T$ is a hyponormal operator\/ $($or\/ if\/ $T$ is a
contraction\/$)$, then
$$
T=I\oplus L
$$
on\/ ${\H=\N(I-T)\oplus\N(I-T)^\perp}\!$ where, in the above decomposition,\/
$I$ is the identity acting on\/ ${\N(I-T)}$ and\/ $L$ is a hyponormal
operator\/ $($or\/ $L$ is a contraction\/$)$ acting on}\/ ${\N(I-T)^\perp\!}$.
\vskip4pt
\item{$\kern-7pt$\rm(b)\kern-1pt}
{\it Moreover, in this case, the power sequence\/ $\{T^n\}$ converges weakly
if and only if the operator\/ $L$ is weakly stable,
$$
T^n\!\wconv P
\quad\hbox{\it if and only if}\quad
L^n\!\wconv O,
$$
so the weak limit\/ ${P\in\BH}$ of\/ $\{T^n\}$, if it exists, is the
orthogonal \hbox{projection}
$$
P=I\oplus O
$$
on\/ ${\H=\N(I-T)\oplus\N(I-T)^\perp}\!$, and hence}\/ $\R(P)=\N({I-T})$.
\end{description}
\vskip-1pt

\vskip5pt
According to (4.2), weak stability of normaloid (and so of hyponormal)
operators implies contractiveness$.$ By Corollary 5.10, if a hyponormal
operator $T$ converges weakly, then it decomposes as ${T\kern-1pt=I\oplus L}$,
where $L$ is a weakly stable hyponormal operator, thus a contraction, and so
is ${T\kern-1pt=I\oplus L}.$ This proves the next result.
\vskip-3pt\noi
$$
\vbox{\noi
{\narrower\narrower
{\it If the power sequence of a hyponormal operator converges weakly, \\
then the hyponormal operator is a contraction}\/.
\vskip0pt}
\vskip-20pt
}                                                                \eqno{(5.6)}
$$
\vskip8pt\noi

%%%%%%%%%%%%%%%%%%%%%%%%%%%  REMARK 5.11
\vskip5pt\noi
{\bf Remark 5.11.}\
Items (ii) and (iii) in Proposition 5.9(a) gave two equivalent conditions for
the weak limit $P$ of $\{T^n\}$ to be orthogonal (i.e., an orthogonal
projection).\ Recall that $P$ always is a projection whenever $T$ is a
bounded operator, according to Proposition 2.6.\ Moreover, it was shown
in \cite[Corollary 3.4]{JJKS} that if the projection $P$ is not orthogonal,
then it is similar to an orthogonal projection:
\vskip2pt\noi
\begin{description}
\item{$\kern-12pt$(a)}
${T^n\kern-2pt\wconv\kern-1ptP}$ if and only if there is an invertible $G$ for
which ${(G^{-1}TG)^n\kern-2pt\wconv\kern-1ptE}$~and \\
${\kern-4ptE}$ is an orthogonal projection.\
$$
\hbox{Since $\;(G^{-1}TG)^n\!=G^{-1}T^nG,\;$ we get $\;E=G^{-1}P\kern1ptG$.}
$$
\end{description}
Corollary 5.10 deals simultaneously with weakly convergent power sequences of
hyponormal operators and of contractions, which boils down to weakly
\hbox{convergent} power sequences of contractions by (5.6).\ The contraction
part of Corollary~5.10~can be extended to operators $T\kern-1pt$ whose power
sequence converges weakly and (instead of ${\|T\|\kern-1pt\le\kern-1pt1})$
${\liminf_n\|T^n\|\kern-1pt\le\kern-1pt1}.$ In fact, the next result was
proved in \cite[\hbox{Theorem 4.12}]{JJKS}.\
\vskip2pt\noi
\begin{description}
\item{$\kern-12pt$(b)}
If ${T^n\!\wconv P}$ and ${\liminf_n\|T^n\|\le1},\;$ then $P$ is
the~orthogonal~projection~onto \\
${\kern-4pt\N(I-T)}=\N(I-T^*)$.\
\end{description}
This is indeed a more general version than its contraction counterpart in
\hbox{Corollary} 5.10.\ The example below showing a weakly stable operator
that is neither uniformly stable nor a contraction with
${\liminf_n\|T^n\|=1}$ was given in \cite[Example~4.14]{JJKS}.
\vskip1pt\noi
\begin{description}
\item{$\kern-12pt$(c)}
There exists a unilateral weighted shift $\kern-1ptT$ with the following
properties.
\begin{description}
\item{$\kern3pt$(i)$\kern3pt$}
${\liminf_n\|T^n\|=1}\;\;$ (so $T\kern-1pt$ is not uniformly stable),
\vskip1.5pt
\item{$\kern1.5pt$(ii)$\kern1.5pt$}
${\limsup_n\|T^n\|=\vartheta}$ for an arbitrary ${\vartheta>1}\;\;$ (so
$T\kern-1pt$ is not a contraction),
\vskip1.5pt
\item{(iii)}
$T$ is weakly stable if and only if it is power bounded.
\end{description}
\vskip1pt\noi
\end{description}
\vskip-1pt

\vskip5pt
Recall from (3.3) that parts of a hyponormal operator are hyponormal.\
Therefore, according to Corollary 5.10 (compare with Corollary 5.2 in light
of (5.6)$\kern.5pt$),
\vskip-3pt\noi
$$
\vbox{\noi
{\narrower\narrower
{\it if a contraction or a hyponormal operator has no identity part, then \\
its power sequence converges weakly if and only if it is weakly stable}\/.
\vskip0pt}
\vskip-20pt
}                                                                \eqno{(5.7)}
$$
\vskip8pt

%%%%%%%%%%%%%%%%%%%%%%%%%%%  REMARK 5.12
\vskip5pt\noi
{\bf Remark 5.12.}\
According to Corollary 5.10, the operator $P$ in Proposition 5.9(a) is an
orthogonal projection if $T$ is a contraction and if $T$ is hyponormal (in
particu\-lar, if $T$ is unitary)$.$ On the other hand, according to
Proposition 5.9, Corollary 5.10 can be further extended to classes of
Hilbert-space operators $T$ for which ${\N(I-T)}$ reduces $T.$ For instance,
Corollary 5.10 can be extended to dominant operators (i.e., Hilbert-space
operators such that ${\R(\alpha I-T)}\sse{\R(\overline\alpha I-T^*)}$ for
\hbox{every} scalar $\alpha$, which defines a class of operators properly
including the hyponormal operators, with the property that ${\N(\alpha I-T)}$
reduces $T$ for every \hbox{scalar ${\alpha\in\FF}$).}

\vskip5pt
Recall from (3.1) that subnormal operators are hyponormal, and from (3.3)
that parts of a subnormal operator are again subnormal$.$ Therefore,
\vskip-3pt\noi
$$
\vbox{\noi
{\narrower\narrower
{\it the characterisations of weak stability for hyponormal operators in
Prop\-osition 5.4(a) and Corollary 5.10 also hold for subnormal operators}\/.
\vskip0pt}
\vskip-20pt
}                                                                \eqno{(5.8)}
$$
\vskip4pt

\vskip5pt\noi
Further conditions for weak stability of subnormal operators were given in
\cite[Prop\-ositions 6.1, 6.4, and Theorem 6.12]{JJKS}$.$ To state them, we
need two \hbox{additional notions}$.$
\vskip5pt\noi
(i) If ${T\kern-1pt\in\BH}$
is subnormal and if $E$ is the $\BK$-valued spectral measure of a mini\-mal
normal extension ${N\in\BK}$ of $T$ (with ${\H\sse\K}$), then the semispectral
measure $F$ of $T$ is defined by ${F(\Lambda)=\varPi E(\Lambda)|_\H}$ in $\BH$
for Borel subsets $\Lambda$ of the complex plane $\CC$, where $\varPi$ is the
orthogonal projection on $\K$ with ${\R(\varPi)=\H}.$ In this case, let 
$F_\TT$ be the restriction of $F$ to the $\sigma$-algebra of Borel sets in the
\hbox{unit circle $\TT.$}
\vskip5pt\noi
(ii) The second necessary notion is that of a Rajchman measure, which
is~any~finite positive (real-valued) measure $\mu$ on the $\sigma$-algebra of
Borel sets in the unit circle~$\TT$~for which ${\int_\TT\!z^k d\mu\to0}$ as
${|k|\to\infty}$ (we will return to Rajchman measures in \hbox{Section~9}).

%%%%%%%%%%%%%%%%%%%%%%%%%%%  PROPOSITION 5.13
\vskip5pt\noi
{\bf Proposition 5.13} \cite{JJKS} (2024).\
{\it Let\/ ${T\kern-1pt\in\BH}$ be a subnormal operator on a Hilbert space
$\H\kern-1pt$ with semispectral measure\/ $F\!$, and let\/ $N$ be a minimal
normal \hbox{extension$\kern-1pt$ of $T\kern-2pt$}\/.}
\begin{description}
\item{$\kern-7pt$\rm(a)}
{\it $T\kern-1pt$ is weakly stable if and only if it is a contraction and\/
${\<F_\TT(\cdot)x\,;x\>}$ is a \hbox{Rajchman} measure for every}\/
${x\in\H}.$
\vskip4pt
\item{$\kern-7pt$\rm(b)}
{\it ${T^n\!\wconv P}$ if and only
if\/ ${T\kern-1pt=I\oplus L}$ on\/ ${\H=\R\oplus\R^\perp\!}$ where\/ $L$ is a
weakly stable subnormal contraction with semispectral measure\/ $G$ such
that\/ ${\<G_\TT(\cdot)v\,;v\>}$ is a Rajchman measure for every\/
${v\in\R^\perp}\!.$ In this case,\/ $P$ is the orthogonal projection on\/ $\H$
with\/ ${\R(P)=\R=\N(I-T)}$}.
\vskip4pt
\item{$\kern-7pt$\rm(c)}
{\it $\{T^n\}$ converges weakly if and only if\/ $\{N^n\}$ converges weakly$.$
In this case,\/ $\N(I-T)=\N(I-N)$ and $N\!=\kern-1pt{I\oplus M}$ on\/
\hbox{$\K\kern-1pt=\kern-1pt\N(I-N)\oplus\N(I-N)^\perp\!$, where}\/ $M$ is a
weakly stable minimal normal extension of the subnormal weakly stable
contraction\/ $L$ of item\/ {\rm(b)}$.$ Also,\/ $T$ is weakly stable if and
only if\/ $N$ is}\/.
\end{description}
\vskip-1pt

\vskip5pt
The decomposition results in \cite[Theorem 3.2, Corollaries 4.4 and 4.5, and
Theorem 6.12]{JJKS} (reviewed above in Propositions 5.8, 5.9, 5.10, and 5.13)
show that those statements referring to operator decomposition also hold, in
addition, if weak convergence/stability is replaced with strong and uniform
convergence/stability.

\vskip5pt
As is well known, uniform stability makes the sequence $\{\|T^n\|\}$ go
exponentially fast to zero$:$ ${T^n\!\uconv O}$ if and only if
${\|T^n\|\le\beta\alpha^n}$ for some ${\beta\ge1}$ and some ${\alpha\in(0,1)}$
(see, e.g., \cite[Proposition 0.4]{MDOT}$\kern.5pt).$ We close this section
with a result from \cite[Corollary 3]{BM} which shows how slow weak (and so
strong) stability may approach zero.

%%%%%%%%%%%%%%%%%%%%%%%%%%%  PROPOSITION 5.14
\vskip5pt\noi
{\bf Proposition 5.14} \cite{BM} (2009).\
{\it $\!$Let\/ ${T\kern-1.5pt\in\kern-.5pt\BH}\kern-.5pt$ on a Hilbert space\/
$\kern-.5pt\H\kern-.5pt$ be weakly stable but not uniformly stable$.$ Let\/
$\kern-.5pt\{\alpha_n\}\kern-.5pt$ be a sequence of positive numbers such
\hbox{that\/ $\kern-.5pt{\alpha_n\!\to\kern-.5pt0}.$} Then for every\/
${\veps\kern-.5pt>\kern-.5pt0}$ there exists an\/
${x\kern-.5pt\in\kern-.5pt\H}$ with\/
\hbox{${0\kern-.5pt<\kern-.5pt\|x\|\kern-.5pt<\kern-.5pt
\sup_n\alpha_n\kern-1.5pt+\kern-.5pt\veps}\kern-.5pt$
for which}}
$$
\alpha_n<|\<T^nx\,;x\>|
\quad\hbox{\it for every}\quad
n\ge0.
$$
\vskip-2pt

%%%%%%%%%%%%%%%%%%%%%%%%%%%%%%%%%%%%%%%%%%%%%%%%%%%%%%%%%  SECTION 6
\section{Ces\`aro Means and Weak Stability}

An operator $T$ acting on a normed space is ergodic if the sequence of
Ces\`aro means $\{\frac{1}{m}\sum_{n=1}^mT^n\}$ converges strongly.\ The Mean
Ergodic Theorem says that {\it \hbox{every} power-bounded operator on a
reflexive Banach space is ergodic} --- see, e.g.,
\cite[Corollary VIII.5.4]{DF}.\ (See also \cite[Theorem 8.22]{EFHN},
\cite[Theorem 2.2]{Yah} and, for a nearly original version of it in terms
of Hilbert-space isometries, \cite[p.16]{Hal0}.)~Thus
$$
T^n\wconv P
\;\limply\;
{\sup}_n\|T^n\|<\infty
\;\limply\;
\smallfrac{1}{m}\hbox{$\sum$}_{n=1}^mT^n\!\sconv E
\;\limply\;
\smallfrac{1}{m}\hbox{$\sum$}_{n=1}^mT^n\!\wconv E
$$
(as ${m\to\infty}$) for operators ${P,E}$ acting on the same reflexive Banach
space as $T\kern-1pt$ acts.

%%%%%%%%%%%%%%%%%%%%%%%%%%%  PROPOSITION 6.1
\vskip5pt\noi
{\bf Proposition 6.1.}\
{\it If\/ $T\kern-1pt$ is an operator acting on a reflexive Banach space\/
$\X\kern-1pt$, then}
\vskip2pt\noi
$$
{\kern-4pt}
\hbox{\it
$T^n\!\!\wconv\!P\;\,\Rightarrow\;\smallfrac{1}{m}
{\sum}_{n=1}^mT^n\!\!\sconv\!E\;$
for projections $\kern-.5ptP,E\kern-.5pt\in\kern-1pt\BX$ with
$\R(E)\kern-1pt=\kern-1pt\R(P),\!\!\!\!\!$}
                                                             \leqno{\rm(a)}
$$
$$
T^n\!\wconv O
\quad\limply\quad
\smallfrac{1}{m}\hbox{$\sum$}_{n=1}^mT^n\!\sconv O.          \leqno{\rm(b)}
$$
{\it If, in addition,\/ $T\kern-1pt$ acts on a Hilbert space, then}
$$
T^n\!\wconv P
\;\;\hbox{\it and the projection $P$ is orthogonal}
\quad\limply\quad
\smallfrac{1}{m}\hbox{$\sum$}_{n=1}^mT^n\!\sconv P.          \leqno{\rm(c)}
$$

\proof
(a)
Take ${T\kern-1pt\in\BX}$ on a reflexive Banach space $\X.$ Suppose
${T^n\!\wconv P}$, consider the Mean~Ergodic Theorem, and let $E$ be the
strong limit of the sequence of Ces\`aro means of $T.$ Then ${E\in\BX}$ (by
the Banach--Steinhaus Theorem) and $E$ is a projection with
${\R(E)=\N({I-T})}$ (see, e.g., \cite[Corollary VIII 5.2]{DF}$\kern.5pt$).\ By
Proposition 2.6(a) and Corollary 2.7, ${P\in\BX}$ is a projection with
${\R(P)=\N(I-T)}.$~So
$$
\R(E)=\R(P)=\N(I-T).
$$
(b)
By the above identity, ${\R(E)=\0}$ if and only if ${\R(P)=\0.}$ Hence
$$
P=O
\quad\iff\quad
E=O.
$$
(c)
By \cite[Corollary VIII 5.2]{DF} we also get ${\R(I-E)=\R(I-T)^-}\!.$
That is,
\vskip2pt\noi
\begin{description}
\item{(i)}
${\N(E)=\R(I-T)^-}$.
\end{description}
\vskip1pt\noi
{From} now on suppose, in particular, that $T$ acts on a Hilbert space.\
Therefore
\vskip1pt\noi
\begin{description}
\item{(ii)}
${\R(I-T)^-\!=\N(I-T^*)^\perp\!}$.
\end{description}
\vskip1pt\noi
If the projection $P$ is orthogonal, then Proposition 5.9(a) ensures that
\vskip1pt\noi
\begin{description}
\item{(iii)}
${\N(I-T^*)=\N(I-T)}$.
\end{description}
\vskip1pt\noi
According to Proposition 2.6(a),
\vskip1pt\noi
\begin{description}
\item{(iv)}
${\N(I-T)=\R(P)}$.
\end{description}
\vskip1pt\noi
Finally, if the projection $P$ is orthogonal, then
\vskip1pt\noi
\begin{description}
\item{(v)}
${\R(P)^\perp\!=\N(P)}$.
\end{description}
By using the above identities sequentially we get
$$
\N(E)=\R(I-T)^-\!=\N(I-T^*)^\perp\!=\N(I-T)^\perp\!=\R(P)^\perp\!=\N(P).
$$
Thus ${\N(E)=\N(P)}$ and ${\R(E)=\R(P)}$ (by item (a)$\kern.5pt).$ This is
enough to ensure that ${E\kern-1pt=\kern-1pt P}$ (because the range and kernel
of projections are algebraic \hbox{complements).\!\qed}

\vskip5pt
The reverse implications in Proposition 6.1 require the notion of subsequence.\

\vskip6pt
Consider the set of all positive integers equipped with its
natural well-ordering, and regard it as a self-indexed sequence of integers so
that ${\NN=\{n\}_{n\ge1}}.$ A subsequence $\{n_k\}=\{n_k\}_{k\ge1}$ of
${\NN=\{n\}_{n\ge1}}$ is precisely a strictly increasing (infinite) sequence
of positive integers$.$ A subsequence $\{a_{n_k}\}$ of an arbitrary $A$-valued
sequence $\{a_n\}$, for an arbitrary nonempty set $A$, is the restriction of
the sequence $\{a_n\}$ to a sub\-sequence $\{n_k\}$ of positive integers
(regarded as a well-ordered \hbox{subset of $\NN$)}.

\vskip5pt
What is behind Proposition 6.1 is the well-known elementary result saying that
if $\X$ is a normed space, $y$ lies in $\X$, and $\{x_n\}$ is an $\X$-valued
sequence, then
$$
x_n\to y
\quad\limply\quad
\smallfrac{1}{m}\hbox{$\sum$}_{n=1}^mx_n\to y
$$
(as ${m\to\infty}$); see, e.g., \cite[Theorem 19.3]{Bar}.\ Thus, if
${T\kern-1pt\in\BX}$ and ${x,y\in\X}$, then
$$
T^nx\to y\!
\quad\limply\quad
\smallfrac{1}{m}\hbox{$\sum$}_{n=1}^mT^nx\to y.
$$
A first Hilbert-space weak version, yielding a result along the same lines as
the above implication, reads as follows.\ Let $\H$ be a Hilbert space, let
$\{x_{n_k}\}$ be a sub\-sequence of an $\H$-valued sequence $\{x_n\}$, and let
$y$ be in $\H.$ Then
$$
x_n\wconv y
\quad\limply\quad
\hbox{there is a subsequence $\{x_{n_k}\}$ such that
$\frac{1}{m}\sum_{k=1}^mx_{n_k}\!\to y$}
$$
(see, e.g., \cite[Exercise 4.25(b)]{Wei}$\kern.5pt).$ Thus, if
${T\kern-1pt\in\BH}$ and ${x,y\in\H}$, then
$$
T^nx\wconv y\!
\quad\limply\quad
\hbox{there is a subsequence $\{T^{n_k}\}$ such that
$\frac{1}{m}\sum_{k=1}^mT^{n_k}x\to y$}.
$$
It was shown in \cite{JK} that the implication becomes an equivalence when the
right-hand side is satisfied for every subsequence $\{T^{n_k}\}$ of $\{T^n\}$,
if $T$ is a contraction.

%%%%%%%%%%%%%%%%%%%%%%%%%%%  PROPOSITION 6.2
\vskip5pt\noi
{\bf Proposition 6.2} \cite{JK} (1971).\
{\it Let\/ $\H$ be a Hilbert space, and let\/ $x$ and\/ $y$ be vectors in\/
$\H.$ If\/ $T\kern-1pt$ is a contraction on\/ $\H$, then}
$$
{\!\!T^nx\wconv y}
\;\,\iff\;\,
\hbox{\it $\frac{1}{m}\sum_{k=1}^mT^{n_k}x\to y$\/
for every subsequence\/ $\{T^{n_k}\}$ of\/ $\{T^n\}$}.             \eqno{(\S)}
$$
\vskip-1pt

\vskip5pt
See also \cite[Theorem II.3.6]{Eis} and the references therein.\
$\kern-1pt$The right-hand side~of $(\S)$ refers to every subsequence of the
power sequence $\{T^n\}$ and not to every subsequence of the sequence of
Ces\`aro means $\{\frac{1}{m}\sum_{k=1}^mT^k\}.$ A similar result for
contractions appeared in \cite[Theorem 1.1]{AS}, apparently independently and
simultaneously, although the convergences in \cite[Theorem 1.1]{AS} are not
shown to be to the same limit $y.\kern-1pt$ The equivalence in $(\S)$,
however, fails if $T$ is just power bounded (rather than a contraction --- see
\cite[Example 2.1]{MT}$\kern.5pt).$ Since the equivalence in $(\S)$ holds
for every $x$ with the same limit $y$ on both sides, we may rewrite
Proposition 6.2~as~follows.

%%%%%%%%%%%%%%%%%%%%%%%%%%%  COROLLARY 6.3
\vskip5pt\noi
{\bf Corollary 6.3.}\
{\it If\/ $T$ is a contraction on a Hilbert space $\H$ and
${P\kern-1pt\in\kern-.5pt\BH}$, then}
$$
T^n\wconv P
\;\,\iff\;\,
\hbox{\it $\frac{1}{m}\sum_{k=1}^mT^{n_k}\sconv P$\/
for every subsequence\/ $\{T^{n_k}\}$\/ of $\{T^n\}$}.
$$
\vskip-4pt

\vskip5pt
As we saw before, the Hilbert-space operator $P$ in Corollary 6.3 is a
projection.\ (Compare Corollary 6.3, which is restricted to contractions,
with Proposition 6.1(a,c).) Recall again that two operators $A$ and $B$ on a
normed space $\X$ are similar if there is an invertible ${G\in\BX}$ (with
${G^{-1}\!\in\BX}\kern.5pt$) such that ${GA=B\kern1pt G}.$

%%%%%%%%%%%%%%%%%%%%%%%%%%%  COROLLARY 6.4
\vskip5pt\noi
{\bf Corollary 6.4.}\
{\it If\/ ${T\kern-1.5pt\in\kern-.5pt\BH}$ is similar to a contraction on
a Hilbert space\/~$\H$,~then}
$$
T^n\wconv P
\;\,\iff\;\,
\hbox{\it $\frac{1}{m}\sum_{k=1}^mT^{n_k}\sconv P\,$
for every subsequence\/ $\{T^{n_k}\}$\/ of $\{T^n\}$},
$$
{\it where\/ ${P\kern-1pt\in\kern-.5pt\BH}$ is a projection}\/.

\proof
Suppose there is a contraction $C$ and an invertible (with a bounded
\hbox{inverse}) $G$ in $\BH$ such that ${G\kern1pt T\kern-1pt=CG}.$ So
${T^n\!\wconv P}$ if and only if ${C^n\!\wconv E}$, with $P={G^{-1}EG}$ and
$E={GPG^{-1}}$ being projections in $\BH$ together.\
(Indeed,~${(G^{-1}CG)^n}=$ ${G^{-1}C^nG\wconv G^{-1}EG=(G^{-1}EG)^2}$; cf.\
Proposition 2.6(b) and Corollary 2.7$.$) Similarly,
$\frac{1}{m}\sum_{k=1}^mT^{n_k}\sconv P\iff
\frac{1}{m}\sum_{k=1}^mC^{n_k}\sconv E.$
Now apply Corollary~6.3.                                                 \qed

\vskip5pt
${\kern-4pt}$Also recall that an operator $T\kern-1pt$ on a normed space $\X$
is power bounded below if there exists $\,{\gamma>0}\,$ such that
$\,{\gamma\|x\|<\|T^nx\|}\,$ for all $\,{n\ge1}\,$ and every
$\,{0\ne x\in\X}.$

\vskip5pt
${\kern-3.5pt}$It was proved in \cite[Theorem 2]{KR} that
\vskip-3pt\noi
$$
\vbox{\noi
{\narrower\narrower
{\it a normed-space operator is power bounded and power bounded below \\
if and only if it is similar to an isometry}\/.
\vskip0pt}
\vskip-20pt
}                                                                \eqno{(6.1)}
$$
\vskip10pt\noi

%%%%%%%%%%%%%%%%%%%%%%%%%%%  COROLLARY 6.5
\vskip5pt\noi
{\bf Corollary 6.5.}\
{\it Let\/ $T\kern-1pt$ be an operator on a Hilbert space\/ $\H$.\ If\/
$T\kern-1pt$ is power bounded and power bounded below, then}
$$
T^n\wconv P
\;\,\iff\;\,
\hbox{\it $\frac{1}{m}\sum_{k=1}^mT^{n_k}\sconv P\,$
for every subsequence\/ $\{T^{n_k}\}$\/ of $\{T^n\}$},
$$
{\it where\/ $P$ is a projection in}\/ $\BH$.

\proof
This is a particular case of Corollary 6.4 according to (6.1). \qed

%%%%%%%%%%%%%%%%%%%%%%%%%%%  REMARK 6.6
\vskip5pt\noi
{\bf Remark 6.6.}\
Under the assumption of Corollary 6.4, viz., if a Hilbert-space opera\-tor
$T\kern-.5pt$ is {\it similar to a contraction}\/, in particular, if it is a
{\it contraction}\/ or if it is {\it similar to an isometry} --- that is, if
it is {\it power bounded and power bounded below}$\kern-1pt$ --- then
$$
T^n\wconv O
\;\;\iff\;\;
\hbox{$\frac{1}{m}\sum_{k=1}^mT^{n_k}\sconv O\,$
for every subsequence $\{T^{n_k}\}$ of $\{T^n\}$},
$$
which, in such particular cases, goes beyond Proposition 6.1(b).

\vskip5pt
It was shown in \cite[Theorem 1.2]{AS} that the contractiveness assumption in
Proposition 6.2 may be replaced by a form of power boundedness below, thus
yielding an apparent generalisation of Proposition 6.5.\ Such replacement and
apparent generalisation, however, are partial because the limit in the next
proposition is not stated to be the same on both sides of the equivalence
$(\ddag)$.

%%%%%%%%%%%%%%%%%%%%%%%%%%%  PROPOSITION 6.7
\vskip5pt\noi
{\bf Proposition 6.7} \cite{AS} (1972).\
{\it If\/ $T$ is an operator on Hilbert space\/ $\H$, and if
$$
\hbox{there exist}\;\;\gamma>0\;\;\hbox{such that}\;\;
\gamma\|x\|<{\limsup}_n\|T^nx\|\;\;\hbox{for every}\;\;0\ne x\in\H,
                                                                \eqno{(\dag)}
$$
\vskip-3pt\noi
then}
\vskip-9pt\noi
$$
\vbox{\hskip-0pt
{\narrower\narrower
{\it $\{T^n\}\;$ converges weakly
$\iff \big\{\frac{1}{m}\sum_{k=1}^mT^{n_k}\big\}\;$ \it converges weakly \\
\phantom{for\/}for every subsequence\/ $\{T^{n_k}\}$ of\/ $\{T^n\}$}\/.
\vskip0pt}
\vskip-20pt
}                                                              \eqno{(\ddag)}
$$
\vskip10pt

\vskip5pt
Observe that power boundedness is implicitly assumed, as it is implied by
weak convergence of $\{T^n\}$, and condition ${(\dag)}$ coincides with
power boundedness below if $T\kern-1pt$ is power bounded, since
${\|T^nx\|\le\|T^{n-m}\|\kern1pt\|T^mx\|}$ for all integers
${0\kern-1pt<\kern-1ptm\le\kern-1ptn}.$

\vskip5pt
Still along the lines of this section, it is worth mentioning an interesting
example exhibited in \cite[Theorem 5.4]{ASu}, where an operator fails to be
strongly stable (in a strongly coercive way) but weak convergence to zero
holds for the Ces\`aro means.

%%%%%%%%%%%%%%%%%%%%%%%%%%%%%%%%%%%%%%%%%%%%%%%%%%%%%%%%%  SECTION 7
\section{Weak Stability for a Class of $\,2\!\times\!2\,$ Operator Matrices}

Let $\K$ and $\Le$ be nonzero complex Hilbert spaces, and consider their
(external) orthogonal direct sum ${\H=\K\oplus\Le}$, again a nonzero complex
Hilbert space$.$ Let $T\kern-1pt$ be an operator on ${\H=\K\oplus\Le}$ given
by the following upper triangular operator~matrix
\vskip5pt\noi
$$
T=\Big(\smallmatrix{A & B \cr
                    O & C \cr}\Big)\!:\K\oplus\Le\to\K\oplus\Le,
$$
\vskip3pt\noi
where ${A\!\!:\K\!\to\K}$, $\,{B\!:\Le\!\to\K}$ and ${C\!:\Le\!\to\Le}$ are
arbitrary bounded linear transforma\-tions$.$ Equivalently, such a model also
applies to operators $T\kern-1pt$ acting on
${\H\kern-1pt=\!\M\oplus\M^\perp}\!$ for some $T$-invariant nontrivial
subspace $\M$ of $\H$ (i.e., for $\M$ being a closed linear manifold of a
Hilbert space $\H$ such that ${\0\ne\M\ne\H}$ and ${T(\M)\sse\M}$), with
${\K\kern-1pt=\kern-1pt\M}$ and ${\Le\kern-1pt=\kern-1pt\M^\perp}\!$, where
${A\kern-1pt=\kern-1ptT|_\M\!:\kern-1pt\M\to\kern-1pt\M}.$ As is readily
verified by~induction,
\vskip5pt\noi
$$
T^n=\Big(\smallmatrix{A^n & B_n \cr
                      O   & C^n \cr}\Big)
\quad\hbox{for every integer}\quad
n\ge0,
$$
\vskip3pt\noi
with ${B_n\!:\Le\!\to\K}$ given by
$$
B_{n+1}=\hbox{$\sum$}^n_{k=0}A^{n-k}B\,C^k
\quad\hbox{for every integer}\quad
n\ge0
\quad\hbox{with}\quad
B_0=O.
$$
\vskip-1pt

%%%%%%%%%%%%%%%%%%%%%%%%%%%  PROPOSITION 7.1
\vskip5pt\noi
{\bf Proposition 7.1.}\
{\it If\/ $T\kern-1pt$ is an operator on\/ ${\H=\K\oplus\Le}$ as above, then
the following assertions are equivalent}\/.
\begin{description}
\item{$\kern-6pt$\rm(a)}
${T^n\!\wconv O}$.
\vskip4pt
\item{$\kern-6pt$\rm(b)}
${A^n\!\wconv O}$, $\;{C^n\!\wconv O},\;$ {\it and}\/ $\,{B_n\!\wconv O}$.
\end{description}

\proof
Take an arbitrary ${x=(u,v)\in\K\oplus\Le}.$ Then, with $B_n$ as given above,
\vskip-2pt\noi
$$
\vbox{
\begin{eqnarray*}
\hskip-20pt
\<T^nx\,;x\>_{_\H}\!
&\kern-6pt=\kern-6pt&
\big\<(A^nu+B_nv,C^nv)\,;(u,v)\big\>_{^\H}\!
=\<A^nu+B_nv\,;u\>_{_\K}\!+\<C^nv\,;v\>_{_\Le}\!                          \\
&\kern-6pt=\kern-6pt&
\<A^nu\,;u\>_{_\K}\!+\<B_nv\,;u\>_{_\K}\!+\<C^nv\,;v\>_{_\Le}.
\end{eqnarray*}
\vskip-19pt
}                                                                  \eqno{(*)}
$$
\vskip7pt\noi
Recall that $\K$, $\Le$ and $\H$ are complex$.$ According to $(*)$, (b)
implies (a) trivially$.$ Con\-versely, by setting ${x=(u,0)}$ and ${x=(0,v)}$
for arbitrary ${u\in\K}$ and ${v\in\Le}$, \hbox{we get}
$$
T^n\wconv O\;\;\limply\;\,A^n\wconv O\;\;\hbox{and}\;\;C^n\wconv O.
$$
Therefore, if (a) holds, then for an arbitrary ${x=(u,v)\in\K\oplus\Le}$,
$$
\<T^nx\,;x\>_{_\H}\!\to0,\,\;\<A^nu\,;u\>_{_\K}\!\to0,
\,\;\hbox{and}\;\;\<C^nv\,;v\>_{_\Le}\!\to0,
$$
and so, by $(*)$, ${\<B_nv\,;u\>_{_\K}\!\to0}$ for every ${u\in\K}$ and
${v\in\Le}.$ Hence, (a) implies (b).                                     \qed

\vskip5pt
Of course, the above argument holds for any sequence of operators converging
to the null operator, not only for power sequences, and not only for
weak~convergence.

\vskip5pt
We discuss below two samples of operators fitting the above pattern whose weak
stability has been considered in the literature.

%%%%%%%%%%%%%%%%%%%%%%%%%%%%%%%%%%%%%%%%%%  SUBSECTION 7.1
\subsection{Brownian-type Operators}

Take an operator
$$
T=\Big(\smallmatrix{V & E \cr
                    O & X \cr}\Big)\,\in\,\B[\K\oplus\Le],
$$
with ${V\in\B[\K]}$, $\,{E\in\B[\Le,\K]}$ and ${X\in\B[\Le]}$ so that,
for every integer ${n\ge0}$,
$$
T^n=\Big(\smallmatrix{V^n & E_n \cr
                      O   & X^n \cr}\Big)\,\in\,\B[\K\oplus\Le],
$$
where, as above, ${E_n\in\B[\Le,\K]}$ is given for every integer ${n\ge0}$ by
$$
E_{n+1}=\hbox{$\sum$}^n_{k=0}V^{n-k}EX^k
\quad\hbox{with}\quad
E_0=O.
$$
According to \cite[Definition 1.1]{CJJS1}, the operator $T$ is of
Brownian type if
\vskip2pt\noi
\begin{description}
\item{$\kern-2pt$\rm(i)\kern6pt}
$V^*V=1\,$
\qquad\qquad\quad
(i.e., $V$ is an isometry),
\vskip4pt
\item{$\kern-4pt$\rm(ii)\kern4pt}
$V^*E=O$
\qquad\qquad\quad
(i.e., $\R(E)\sse\N(V^*)$; equivalently, $\R(V)\perp\R(E)$),
\vskip4pt
\item{$\kern-6pt$\rm(iii)\kern2pt}
$XE^*E=E^*EX$
\qquad
(i.e., $X$ commutes with the nonnegative $E^*E$).
\end{description}
\vskip-1pt

\vskip5pt
$\!$Weak stability of Brownian-type operators was investigated in
\cite[Corollary~2.2]{CJJS2}.

%%%%%%%%%%%%%%%%%%%%%%%%%%%  PROPOSITION 7.2
\vskip5pt\noi
{\bf Proposition 7.2} \cite{CJJS2} (2024).\
{\it $\!$A Brownian-type operator $T\kern-1.5pt$ is weakly stable if
\hbox{and $\kern-1pt$only} if\/ $X$ and the unitary part of the isometry\/
$V\!$ are weakly stable and\/ $\{E_n\}$ is
bounded}\/.

\vskip5pt
Suppose (i), (ii), (iii) hold true$.$ The proof of Proposition 7.2 in
\cite[Corollary 2.2]{CJJS2} uses Proposition 7.1 with $E_n$ given above, from
which it can be easily verified~that
$$
E_{n+1}^*E_{n+1}
=\hbox{$\sum$}_{k=0}^nX^{*k}X^kE^*E
=E^*E\,\hbox{$\sum$}_{k=0}^nX^{*k}X^k
$$
for every nonnegative integer $n$, and it also uses the equivalence
$$
E_n\!\wconv O
\quad\iff\quad
{\sup}_n\|E_n\|<\infty
$$
proved in \cite[Lemma 3.4]{CJJS2}$.$ Proposition 7.2 can be rewritten as
follows.

%%%%%%%%%%%%%%%%%%%%%%%%%%%  COROLLARY 7.3
\vskip5pt\noi
{\bf Corollary 7.3.}\
{\it A Brownian-type operator\/ $T$ is weakly stable if and only if\/ $X$ and
$V\!$ are weakly stable and\/ ${\sum_{k=0}^\infty\|EX^kv\|^2\!<\infty}$ for
every}\/ ${v\in\Le}$.
\vskip-2pt\noi

\proof
As we saw in (4.7), an isometry (as is $V$ in (i)$\kern.5pt$) is weakly stable
if and only if its unitary part is weakly stable$.$ Moreover,
\begin{eqnarray*}
\|E_{n+1}v\|^2
&\kern-6pt=\kern-6pt&
\<E_{n+1}^*E_{n+1}v\,;v\>
=\!\hbox{$\sum$}_{k=0}^n\<X^{*k}X^kE^*Ev\,;v\>
=\!\hbox{$\sum$}_{k=0}^n\<X^kE^*Ev\,;X^kv\>                                \\
&\kern-6pt=\kern-6pt&
\hbox{$\sum$}_{k=0}^n\<E^*EX^kv\,;X^kv\>
=\hbox{$\sum$}_{k=0}^n\<EX^kv\,;EX^kv\>
=\hbox{$\sum$}_{k=0}^n\|EX^kv\|^2
\end{eqnarray*}
for each ${n\ge0}$ and every ${v\in\Le}$ and so, by the Banach--Steinhaus
Theorem,
\goodbreak\noi
\begin{eqnarray*}
{\sup}_n\|E_n\|<\infty
&\kern-6pt\iff\kern-6pt&
{\sup}_n\|E_nv\|<\infty\quad\forall v\in\Le                                \\
&\kern-6pt\iff\kern-6pt&
{\sup}_n\hbox{$\sum$}_{k=0}^n\|EX^kv\|^2<\infty\quad\forall v\in\Le        \\
&\kern-6pt\iff\!\kern-6pt&
\hbox{$\sum$}_{k=0}^\infty\|EX^kv\|^2<\infty\quad\forall v\in\Le. 
\end{eqnarray*}
\vskip-4pt\noi
Now apply Proposition 7.2.                                               \qed

%%%%%%%%%%%%%%%%%%%%%%%%%%%  COROLLARY 7.4
\vskip5pt\noi
{\bf Corollary 7.4.}\
$\kern-.5pt${\it Let\/ $T\kern-1.5pt$ be a Brownian-type operator.\
$\kern-.5pt$Suppose\/ $\kern-.5pt X\!$ is a quasinormal con\-traction,
${E\kern-.5pt=\kern-.5pt(I\kern-1pt-\!X^*\!X)^{1/2}}\!$, and $V\!$ is an
isometry such that ${(I\kern-1pt-\!X^*\!X)\,V\kern-2.5pt=\kern-.5ptO}.\!$~Then
\vskip5pt\noi
\centerline{$T$ is weakly stable $\iff$ $V$ and\/ $X$ are weakly stable\/.}}
\vskip-2pt\noi

\proof
Since $X$ is a contraction if and only if ${I-X^*\!X}$ is a nonnegative
operator, set ${E=(I-X^*\!X)^{1/2}}\kern-1pt.$ Since ${(I-X^*\!X)\,V\!=O}$ if
and only if ${E^*V\!=O}$ (because $\N(E^2)=\N(E)$), assumption (ii) is
satisfied$.$ Since $X$ is quasinormal if and only if
${(I-X^*\!X)X}={X(I-X^*\!X)}$, assumption (iii) is satisfied$.$ Also, for
every ${u\in\K}$,
$$
\|Eu\|^2=\<(I-X^*X)u\,;u\>=\|u\|^2-\|Xu\|^2.
$$
As $X$ is a contraction, there is a nonnegative contraction $A$ for which
${X^{*n}X^n\sconv A}$ (see, e.g., \cite[Proposition 3.1]{MDOT}$\kern.5pt$) so
that ${\|X^nu\|\to\|A^{1/2}u\|}$ for every ${u\in\K}.$ Hence
\begin{eqnarray*}
\hbox{$\sum$}^\infty_{k=0}\|EX^ku\|^2
&\kern-6pt=\kern-6pt&
\hbox{$\sum$}^\infty_{k=0}\big(\|X^ku\|^2-\|X^{k+1}u\|^2\big)              \\
&\kern-6pt=\kern-6pt&
\|u\|^2-{\lim}_n\|X^{n+1}u\|^2
=\|u\|^2-\|A^{1/2}u\|^2\le\|u\|^2,
\end{eqnarray*}
and therefore ${\sum}^\infty_{k=0}\|EX^ku\|^2<\infty$, for every $u\in\K.$ Now
apply Corollary 7.3.                                                      \qed

\vskip5pt
The special case of $X$ quasinormal was systematically considered in
\cite{CJJS1}.

%%%%%%%%%%%%%%%%%%%%%%%%%%%%%%%%%%%%%%%%%%  SUBSECTION 7.2
\subsection{Weak Quasistability}

To proceed, we need the notion of weak quasistability$.$

\vskip6pt
Recall that an operator $T\kern-1pt$ on a normed space
$\kern-.5pt\X\kern-.5pt$ is weakly stable if, for every ${x\in\X}$,
$\,\lim_n\kern-1pt|f(T^nx)|\kern-1pt=\kern-1pt0$ for every ${f\in\X^*}\!$
(i.e., ${w\hbox{\,-}\lim_nT^nx=0}$ for every ${x\in\X}).\!$
\begin{description}
\item{$\kern-7pt$\small(1)}
A normed-space operator ${T\kern-1pt\in\BX}$ is {\it weakly quasistable}\/ if
$$
\hbox{for every}\;\,
x\in\X,
\quad
{\liminf}_n|f(T^nx)|=0
\quad\hbox{for every}\;\,
f\in\X^*
$$
(i.e.,
${w\hbox{\kern1pt-\kern-.5pt}\liminf_n\kern-.5pt T^nx\kern-.5pt=\kern-.5pt0}$
for every ${x\kern-.5pt\in\kern-.5pt\H}).$ In particular, in a Hilbert
space~$\H$,
$$
{\liminf}_n|\<T^nx\,;y\>|=0
\quad\hbox{for every}\quad
x,y\in\H.
$$
\end{description}
\vskip-4pt\noi
Equivalently,
\begin{description}
\item{$\kern-7pt$\small(2)}
an operator $T\kern-1pt$ is weakly quasistable if, for every $x$, there is a
subsequence $\{T^{n_k}\}$ of $\{T^n\}$ (that may depend on $x$) such that
${w\hbox{\,-}\lim_kT^{n_k}x=0}.$
\end{description}
In other words,
\begin{description}
\item{$\kern-7pt$\small(3)}
${T\kern-1pt\in\BX}$ is weakly quasistable if there is a subsequence
$\{T^{n_k}\}$ of $\{T^n\}$ (that may depend on ${x\in\X}$ but not on
${f\in\X^*}$) such that
$$
\hbox{for every}\;\,
x\in\X,
\quad
{\lim}_k|f(T^{n_k}x)|=0
\quad\hbox{for every}\;\,
f\in\X^*\!.
$$
(i.e.,
${w\hbox{\kern1pt-\kern-.5pt}\lim_k\kern-.5pt T^{n_k}x=0}$ for every
${x\kern-.5pt\in\kern-.5pt\H}).$ Again, in a Hilbert space $\H$,
$$
{\lim}_k|\<T^{n_k}x\,;y\>|=0
\quad\hbox{for every}\quad
x,y\in\H.
$$
\end{description}
If there exists such a subsequence $\{T^{n_k}\}$ of $\{T^n\}$ for a given $x$,
then $\{T^{n_k}\}$ is referred to as a {\it subsequence of weak
stability}\/ for $x.$ Thus,
\begin{description}
\item{$\kern-7pt$\small(4)}
an operator $T\kern-1pt$ is weakly quasistable if and only if there is at
least one subsequence $\{T^{n_k}\}=\{T^{n_k(x)}\}$ of weak stability for each
$x$.
\end{description}
Therefore, a subsequence $\{T^{n_k}\}$ of $\{T^n\}$ is not of weak
stability for ${x\in\X}$ (or ${x\in\H}$) if there exists ${f\in\X^*}$
(or ${y\in\H}$) such that
$$
{\lim}_k|f(T^{n_k}x)|\ne0
\qquad
(\,{\rm or}\;\,{\lim}_k|\<T^{n_k}x\,;y\>|\ne0\,).
$$
In such a case, the operator $T\kern-1pt$ is said to have a {\it subsequence\/
$\{T^{n_k}\}\kern-1pt$ of weak instability}\/ for ${x\in\X}$ (or ${x\in\H}).$
Equivalently,
\begin{description}
\item{$\kern-7pt$\small(5)}
$T\kern-1pt$ has a subsequence of weak instability for ${x\in\X}$ (or
${x\in\H}$) if and only if there is a subsequence $\{T^{n_k}\}$ of $\{T^n\}$
and ${f\in\X^*}$ (or ${y\in\H}$) such~that
$$
{\limsup}_k|f(T^{n_k}x)|>0
\qquad
(\,{\rm or}\;\,{\limsup}_k|\<T^{n_k}x\,;y\>|>0\,).
$$
\end{description}
This, in turn, is equivalent to saying that
\begin{description}
\item{$\kern-7pt$\small(6)}
there exists ${f\in\X^*}$ (or ${y\in\H}$) and a subsequence
${\{T^{n_i}\}=\{T^{{n_k}_i}\}}$ of a subsequence $\{T^{n_k}\}$ of weak
instability for ${x\in\X}$ (or ${x\in\H}$) such that
$$
{\liminf}_i|f(T^{n_i}x)|>0
\qquad
(\,{\rm or}\;\,{\liminf}_i|\<T^{n_i}x\,;y\>|>0\,).
$$
\end{description}
Such a subsequence $\{T^{n_i}\}\kern-1pt$, which is a subsequence of a
subsequence of weak instability, is referred to as a {\it subsequence of
strict weak instability}\/ for ${x\in\X}$ (or~${x\in\H})$.

\vskip5pt
Weak stability trivially implies weak quasistability.\ The converse fails
even for \hbox{Hilbert}-space power-bounded operators, as we will see in
Proposition 7.5.\ \hbox{Differently} from weak stability, weak quasistability
does not imply power boundedness \cite[Remark 4.4]{Kub5}.\ Also, the
counterparts of uniform and power-bounded strong quasi\-stabilities coincide
with uniform and strong stabilities \cite[Propositions~4.1,~4.2]{KV4}.

%%%%%%%%%%%%%%%%%%%%%%%%%%%%%%%%%%%%%%%%%%  SUBSECTION 7.3
\subsection{The Foguel Operator}

Let $S_+\!$ be the unilateral shift of multiplicity one shifting an
orthonormal basis $\{e_k\}$ for an infinite-dimensional separable Hilbert
space $\K$ and take its adjoint ${S_+}^{\!\!\!*}.$ Let $P$ on $\K$ be an
orthogonal projection onto $\R(P)$ $=\span\{e_j\!:j\in\JJ\}^-\!$, where $\JJ$
is a sparse infinite set of positive integers such that if ${i,j\in\JJ}$ with
${i<j}$, then ${2i<j}$ (e.g., $\JJ={\{j\ge 1\!:\;j=3^k;\;k\ge 0\}}$, the set
of all integral powers of 3)$.$ The Foguel operator on the Hilbert space
${\H=\K\oplus\K}$~\cite{Fog2,Hal1},
$$
F=\Big(\smallmatrix{{S_+}^{\!\!\!*} & P \cr
                    O & S_+ \cr}\Big)\!:\H=\K\oplus\K\to\H=\K\oplus\K,
$$
\vskip1pt\noi
was the first example of a power-bounded operator not similar to a
contraction.

\vskip5pt
The power-bounded operator $F$ was built to satisfy
${\Z(F)\cap\Z(F^*)^\perp\!\ne\kern-1pt\0}$, with $\Z(F)$ as in
Proposition 5.1, so that $F$ is not weakly stable (because
${\Z(F^*)^\perp\!\ne\0}$ and $F\kern-1pt$ and $F^*\!$ are weakly stable
together --- see, e.g.,
\cite[Proposition 8.1 and Remark 8.7]{MDOT}$\kern.5pt).$ But $F$ is weakly
quasistable as shown in \cite[Proposition~4.3]{KV4}.

%%%%%%%%%%%%%%%%%%%%%%%%%%%  PROPOSITION 7.5
\vskip5pt\noi
{\bf Proposition 7.5} \cite{KV4} (2024).\
{\it The power-bounded and weakly unstable Foguel operator $F$ is weakly
quasistable\/.}

%%%%%%%%%%%%%%%%%%%%%%%%%%% REMARK 7.6
\vskip5pt\noi
{\bf Remark 7.6.}\
Besides being power bounded (${\sup_n\|F^n\|<\infty}$), not similar to any
contraction (${\|GFG^{-1}\|>1}$ for every invertible operator $G$ with a
bounded inverse $G^{-1}$), weakly unstable (${w\hbox{\,-}\lim_nF^nx\ne0}$ for
some ${x\in\H}$), and weakly quasistable (${w\hbox{\,-}\liminf_nF^nx=0}$ for
every ${x\in\H}$), the Foguel operator $F$ is such that
$$
1=r(F)<w(F)<\|F\|,                                            \leqno{\rm(a)}
$$
where ${w(T)=\sup_{\alpha\in W(T)}|\alpha|}$ is the numerical radius of an
operator ${T\kern-1pt\in\BH}$, with
${W(T)\kern-1pt=\kern-1pt
\{\alpha\kern-1pt\in\FF\!:\alpha\kern-1pt=\kern-1pt\<Tx;x\>
\;\hbox{for some $x\in\kern-1pt\H$ with $\|x\|\kern-1pt=\kern-1pt1$}\}}$
standing for the numerical range of $T.$ As is well known, the numerical
radius $w(T)=\sup_{\|x\|=1}|\<Tx\,;x\>|$ is such that
${r(T)\le w(T)\le\|T\|\le2w(T)}$ for \hbox{every} ${T\kern-1pt\in\BH}.$ (For
treatises on numerical range, see, e.g., \cite{GR,WG}$.$) To confirm the
expression in (a), observe that the identity holds because $F$ is
power-bounded and uniformly unstable, the first strict inequality comes from
\cite[Theorem 3.4]{GLWW} (since $P$ is an orthogonal projection --- and so
${\|P\|=1}$ --- which can be viewed as a noncompact diagonal operator), and
the second strict inequality comes from a well-known result ensuring that
${w(T)=\|T\|}$ implies ${r(T)=w(T)}$ for every
${T\kern-1pt\in\BH}$ (see, e.g., \cite[Theorem 2.16]{ST2}$\kern.5pt).$ As a
matter of fact, it was shown in \cite[Proposition 1.2]{GLWW} that
${w(F)\le\smallfrac{3}{2}}$ and in \cite[Example 2.6]{Gar} that
$\|F\|=\frac{1+\sqrt5}{2}.$ Thus (a) can be strengthened to
$$
r(F)=1<w(F)\le\smallfrac{3}{2}
=\smallfrac{1+\sqrt4}{2}<\smallfrac{1+\sqrt5}{2}=\|F\|.
$$
Moreover, the Foguel operator also contributes to the line of thought in
Remark 5.11(c).\ In fact, it is a power bounded noncontraction (actually, not
similar to any contraction), but this time a weakly unstable operator, such
that
$$
1={\liminf}_n\|F^n\|<{\limsup}_n\|F^n\|,                      \leqno{\rm(b)}
$$
the proof of which is easy by applying the argument in
\cite[proof of Proposition~4.3]{KV4}.

\vskip5pt
Proposition 7.5 says that, although there is a subsequence $\{F^{n_i}\}$ of
$\{F^n\}$ such that ${\<F^{n_i}x\,;y\>\not\to0}$ for some ${x,y\in\H}$, it is
true that for each ${x\in\H}$ there is a subsequence $\{F^{n_k}\}$ of
$\{F^n\}$ (depending on $x$) such that ${\<F^{n_k}x\,;y\>\to0}$ for
every~${y\in\H}$.

\vskip5pt
Since the unilateral shift is weakly stable (and so is its adjoint --- cf$.$
(3.4)$\kern.5pt$) and $F$ is weakly unstable, it follows by Proposition 7.1
that $F$ fails to be weakly stable be\-cause the sequence $\{P_n\}$, with
${P_{n+1}\!=\!\sum_{k=0}^n{S_+}^{\!\!\!\!*\,n-k}P{S_+}^{\!\!\!\!k}}$ for
${n\kern-1pt>\kern-1pt0}$, does not converge weakly to the null operator$.$
And the reason for such a non-convergence is the fact that $\JJ$ is too sparse
a set, so that every subsequence $\{n_j\}$ of the positive inte\-gers indexed
by $\JJ$ is such that the increments ${n_{j+1}-n_j}$ increase unboundedly as
${j\to\infty}.$ This property was fundamental for proving that $F$ is not
weakly stable (see, e.g., \cite[Remark 8.7]{MDOT}$\kern.5pt).$ In other words,
the lack of {\it boundedly spaced}\/ subsequences $\{F^{n_j}\}$ indexed by
$\JJ$ was a reason for the weakly quasistable $F$ to be weakly unstable.

%%%%%%%%%%%%%%%%%%%%%%%%%%%%%%%%%%%%%%%%%%%%%%%%%%%%%%%%%  SECTION 8
\section{Boundedly Spaced Subsequences and Weak Quasistability}

A subsequence $\{n_k\}$ of the sequence of positive integers is of
{\it bounded increments}\/ (or of {\it bounded gaps}\/) if
${\sup_k(n_{k+1}-n_k)<\infty}.$ A subsequence $\{a_{n_k}\}$ of an arbitrary
$A$-valued sequence $\{a_n\}$ (for an arbitrary nonempty set $A$) is
{\it boundedly spaced}\/ if it is indexed by a subsequence $\{n_k\}$ of
bounded increments.

\vskip5pt
The next result from \cite[Theorem 5.3]{KV4} gives a sufficient condition for
weak stability in a normed space in terms of weak quasistability (as defined
in \hbox{Subsection~7.2}).

%%%%%%%%%%%%%%%%%%%%%%%%%%%  PROPOSITION 8.1
\vskip5pt\noi
{\bf Proposition 8.1} \cite{KV4} (2024).\
{\it If a normed-space operator\/ ${T\kern-1pt\in\kern-.5pt\BX}$ has a
boundedly spaced subsequence of weak quasistability for every\/
${x\in\kern-1pt\X}$, then\/ $T\kern-1pt$ is weakly stable}\/.

\vskip5pt
The tensor product counterpart of Remark 4.1 is not as neat as the direct sum:
it requires the notion of boundedly spaced subsequences of power sequences$.$
Let~${\K\kern-1pt\hotimes\Le}$ denote the completion of the tensor product of
Hilbert spaces $\K$ and $\Le$ equipped with the natural reasonable crossnorm
induced by inner products, and let ${A\hotimes B}$ be the extension over
completion of the tensor product of operators $A$ in $\BK$ and $B$ in $\BL$
(see, e.g., \cite[Chapter 9]{BMTP}$\kern.5pt$).\ According to
\cite[Theorems 2 and 3]{KV2}, weak sta\-bility is transferred from operators
to their tensor product quite nicely, but the con\-verse requires a boundedly
spaced condition.\ Indeed, as is well known (see, e.g.,
\cite[Definition 3.16, Theorems 3.3 and 3.4]{BMTP}$\kern.5pt$), for every
${x,u\kern-1pt\in\K}$ and ${y,v\kern-1pt\in\Le}$,~and~${n\kern-.5pt\in\NN}$,
$$
|\<(A\hotimes B)^n(x\otimes y)\,;u\otimes v\>|=|\<A^nx\,;u\>|\,|\<B^ny\,;v\>|,
$$
$$
\|(A\hotimes B)^n\|=\|A^n\|\,\|B^n\|.
$$
\vskip3pt\noi
This is enough to ensure that if one of $A$ or $B$ is weakly stable (thus
power bounded) and the other is power bounded, then ${A\hotimes B}$ is weakly
stable$.$ Conversely, suppose ${A\hotimes B}$ is weakly stable (thus power
bounded) and $A$ and $B$ are power bounded$.$~If $A$ has a boundedly spaced
subsequence of strict weak instability, then there~exist $u_0$ in $\K$ and a
boundedly spaced subsequence $\{A^{n_i}\}$ of weak instability for some $x_0$
in $\H$ such that
$0<\liminf_i|\<A^{n_i}x_0\,;u_0\>|\le\limsup_i|\<A^{n_i}x_0\,;u_0\>|<\infty.$
%%%%%%%%%%%%%%%%%%%%%%%%%%%  FOOTNOTE
$\kern-3pt$\footnote{$\kern2pt$There is a typo in \cite[p.410, line -6]{KV2}
and \cite[p.233, line 24]{BMTP}:\ $\limsup_k$ should read $\liminf_k$.}
%%%%%%%%%%%%%%%%%%%%%%%%%%%
Therefore, since ${\lim_k|\<A^{n_i}x_0\,;u_0\>|}\,{|\<B^{n_i}y\,;v\>|}=0$ for
every ${y,v}$ in $\Le$ (as ${A\hotimes B}$ is weakly~stable), we get
${\lim_i|\<B^{n_i}y\,;v\>|=0}$ for every ${y,v}$ in $\Le$, which implies that
$\{B^{n_i}\}$ is a boundedly spaced subsequence (because $\{n_i\}$ is of
bounded increments) of weak quasistability for every ${y\in\Le}.$ So $B$
is weakly stable by Proposition 8.1$.\kern-1pt$ This is essentially the
argument in the proof of \cite[Theorem 3]{KV2}$.$ Both results are stated
together~\hbox{below}.

%%%%%%%%%%%%%%%%%%%%%%%%%%%  PROPOSITION 8.2
\vskip6pt\noi
{\bf Proposition 8.2} \cite{KV2} (2008).\
{\it Let\/ $A$ and\/ $B$ be Hilbert-space operators on\/ $\K$ and\/ $\Le$,
respectively, and consider their tensor product\/ ${A\hotimes B}$ on}\/
${\K\hotimes\Le}$.
\vskip2pt\noi
\begin{description}
\item{$\kern-9pt$\rm(a)}
{\it If one of\/ $A$ or\/ $B$ is weakly stable and the other is power bounded,
then their tensor product\/ ${A\hotimes B}$ is weakly stable}\/.
\vskip4pt
\item{$\kern-9pt$\rm(b)}
{\it Conversely, suppose\/ ${A\hotimes B}$ is weakly stable and both\/ $A$
and\/ $B$ are power bounded.\ If\/ one of\/ $A$ or\/ $B$ is not weakly stable
but has at least one boundedly spaced subsequence of strict weak instability,
then the other is weakly stable}\/.
\end{description}
\vskip-1pt

\vskip5pt
A notion similar to weak quasistability was considered in
\cite[Definition 1.2]{ES}, namely, a power bounded Hilbert-space operator
${T\kern-1pt\in\BH}$ is {\it almost weakly stable}\/~if
$$
\hbox{$0$ is a weak accumulation point of every orbit under $T$}.
$$
\vskip-3pt

%%%%%%%%%%%%%%%%%%%%%%%%%%%  REMARK 8.3
\vskip5pt\noi
{\bf Remark 8.3.}\
Let $\X$ be a normed space$.$ The orbit of a vector ${y\in\kern-1pt\X}$ under
an operator ${T\kern-1pt\in\BX}$ is the set
$$
\Oe_T(y)
=\{{T^ny\in\X}\!:\hbox{for every integer $n\kern-1pt\ge\kern-1pt0$}\}.
$$
Let $\Oe_T(y)^{-w}\!$ be the weak \hbox{closure} of $\Oe_T(y)$ (i.e., the
closure of $\Oe_T(y)$ in the weak topology on $\X).\kern-1pt$ So $0$ is a weak
accumulation point of \hbox{every} $\Oe_T(y)$ if
${0\kern-1pt\in\kern-1pt\Oe_T(y)^{-w}}\kern-1pt$~for \hbox{every} ${y\in\X}.$
Since the weak topology is not metrisable, a notion different from weak
closure is that of weak limit set.\ The {\it weak limit set}\/
$\Oe(y)^{-wl}\!$ of $\Oe_T(y)$ is the set of all weak limits of weakly
convergent $\Oe_T(y)$-valued sequences;~that is,
$$
\Oe_T(y)^{-wl}\!=
\{{x\in\kern-1pt\X\!:x\kern-1pt=\kern-1ptw\hbox{\,-}{\lim}_kx_k\;
\hbox{with}\;x_k\!\in\Oe_T(y)}\}
=\{{x\in\kern-1pt\X\!:x\kern-1pt=\kern-1ptw\hbox{\,-}{\lim}_kT^{n_k}y}\}.
$$
As we saw in subsection 7.2, weak quasistability for an operator $T\kern-1pt$
means that, for every ${y\in\kern-1pt\X}\kern-1pt$, there~is~a sub\-sequence
$\{T^{n_k}\}$ of $\{T^n\}$ such that
${w\hbox{\,-}\lim_kT^{n_k}y\kern-1pt=\kern-1pt0}.$ So
$$
\hbox{$T\kern-1pt$ is weakly quasistable if and only if
${0\kern-1pt\in\kern-1pt\Oe_T(y)^{-wl}}$ for every ${y\in\X}.$}
$$
But ${\Oe_T(y)^{-wl}\!\sse\kern-1pt\Oe_T(y)^{-w}}$; the inclusion may be
proper \cite[Proposition 2.1]{Kub4}.\ Thus,
$$
{\kern-80pt}
\hbox{\it a weakly quasistable operator is almost weakly stable}\/.
                                                                 \eqno{(8.1)}
$$

\vskip1pt
The results in Proposition 8.4 below on almost weak stability are from
\cite[Prop\-ositions 2.3, 3.4, 4.1 and Theorems 2.4, 3.5, 4.3]{ES}.

%%%%%%%%%%%%%%%%%%%%%%%%%%%%%% PROPOSITION 8.4
\vskip5pt\noi
{\bf Proposition 8.4} \cite{ES} (2008).\
{\it Consider the definition of almost weakly stable operators in $\BH$,
acting on a separable Hilbert space}\/ $\H$.
\vskip2pt\noi
\begin{description}
\item{$\kern-13pt$\rm(${\scriptstyle\U}_1$)\kern-2pt}
{\it The set\/ $\W_U\!$ of all almost weakly stable unitary operators is
residual $($dense, {\rm as a countable intersection of open and dense sets;
i.e., it includes a dense\/~$G_\delta$)} in the set\/ $\U\kern-.5pt$ of all
unitary operators equipped with the strong* operator topology {\rm (which
is complete under an appropriate metric)}\/.}
\vskip4pt\noi
(The strong* $\!$operator topology is defined here as the topology
generated~by the family of seminorms $\{p_x\}_{x\in\H}$ where
${p_x(T)\kern-1pt=\kern-1pt(\|Tx\|^2\!+\kern-1pt\|T^*x\|^2)^{1/2}}\!.$
Convergence of $\{T_n\}$ to $T$ in this topology means
${\|(T_n-T)x\|^2}\!+{\|(T_n-T)^*x\|^2}\!$ $\to0$ for every ${x\in\H}$;
equivalently, ${T_n\!\sconv T}$ and ${T_n^*\!\sconv T^*}\!.$ Recall:\ if an
operator $T$ is normal (in particular, unitary), then ${\|Tx\|=\|T^*x\|}$
for every ${x\in\H}$.)
\vskip4pt
\item{$\kern-13pt$\rm(${\scriptstyle\U}_2$)\kern-2pt}
{\it The set\/ $\Se_U\!$ of all weakly stable unitary operators is of first
category in\/ $\U$}\/.
\vskip4pt
\item{$\kern-13pt$\rm(${\scriptstyle\Ie}_1$)\kern-2pt}
{\it The set\/ $\W_I\kern-1pt$ of all almost weakly stable isometries is
residual $($dense$)$ in the set\/ $\Ie$ of all isometries equipped with the
strong operator~topology {\rm (which is complete under an appropriate
metric)}\/.}
\vskip4pt
\item{$\kern-13pt$\rm(${\scriptstyle\Ie}_2$)\kern-2pt}
{\it The set\/ $S_I\kern-1pt$ of all weakly stable isometries is of first
category in\/ $\Ie$}\/.
\vskip4pt
\item{$\kern-13pt$\rm(${\scriptstyle\C}_1$)\kern-2pt}
{\it The set $\W_C\kern-1pt$ of all almost weakly stable contractions is
residual in the set $\C$ of all contractions equipped with the weak operator
topology {\rm (which is complete under an appropriate metric)}\/.}
\vskip4pt
\item{$\kern-13pt$\rm(${\scriptstyle\C}_2$)\kern-2pt}
{\it The set\/ $\Se_C\kern-1pt$ of all weakly stable contractions is of first
category in\/ $\C$}\/.
\vskip4pt
\item{$\kern-13pt$\rm(${\scriptstyle\C}_3$)\kern-2pt}
{\it The set\/ $\W_U\!$ of all almost weakly stable unitary operators is dense
in}\/ $\C$.
\end{description}

%%%%%%%%%%%%%%%%%%%%%%%%%%%  REMARK 8.5
\vskip5pt\noi
{\bf Remark 8.5.}\
We close this section by remarking on some additional results on weak
stability and its weaker forms beyond those considered above$.$ Let $\X$ be a
normed space$.$ The projective orbit of a vector ${y\in\kern-1pt\X}$ under an
operator ${T\kern-1pt\in\BX}$ is the orbit of span of $\{y\}$; that is, the
orbit $\Oe_T([y])$ of $[y]$:
$$
\Oe_T([y])
=\Oe_T(\span\{y\})
=\big\{\alpha T^ny\in\X\!:\,\alpha\in\FF,\;n\in\NN_0\big\}.
$$
An operator $T$ is a {\it weakly l-sequentially supercyclic operator}\/ if
$$
{\Oe_T([y])^{-wl}\kern-1pt=\kern-1pt\X}
$$
for some ${y\in\kern-1pt\X}$ (i.e., the weak limit set of $\Oe_T([y])$ is the
whole space $\X$; cf.\ Remark 8.3).\ Equivalently, if there exists ${y\in\X}$
such that for every ${x\in\X}$ there is a scalar sequence $\{\alpha_k\}$ such
that ${\alpha_kT^{n_k}y\wconv x}$ for some subsequence $\{T^{n_k}\}$ of
$\{T^n\}.$ That {\it every weakly l-sequentially supercyclic power bounded
operator on a normed space is weakly quasistable}\/ was shown in
\cite[Corollary 4.3]{KD} (and it was asked in \cite{KV4} {\it whether the
Foguel operator of\/ Proposition 7.5 is weakly l-sequentially
\hbox{supercyclic}}\/)$.$~For~ap\-pli\-cations of boundedly spaced
subsequences and weak quasistability towards weak l-se\-quential
supercyclicity and weak stability for unitary and power bounded oper\-ators,
see \cite[\hbox{Theorem 2}]{KV3} and \cite[Theorem 6.2]{KV4}.\ For a result
equivalent to Propo\-sition 8.1 on a Banach space, see
\cite[Theorem II.3.4]{Eis}$.\kern-1pt$ Sufficient conditions for weak
stability, based on resolvent operators, were presented in
\cite[Theorem II.3.13]{Eis}$.$~Fur\-ther formulations of weaker forms of weak
stability (such as almost weak stability) were considered in
\cite[Theorem II.4.1, Remark II.4.4, Corollary~II.4.9]{Eis}.

%%%%%%%%%%%%%%%%%%%%%%%%%%%%%%%%%%%%%%%%%%%%%%%%%%%%%%%%%  SECTION 9
\section{Weak Stability of Unitary Operators}

Weakly stable normaloid operators are contractions (4.2), and Hilbert-space
contractions are weakly stable if and only if their unitary parts are weakly
\hbox{stable} (5.1).\ So weak stability of Hilbert-space operators in these
classes is reduced to~weak~sta\-bility of unitary operators, and Foguel
decomposition (Proposition 5.1) places uni\-tary operators at the centre
of the weak stability discussion.\ This section highlights the major role
played by unitary operators in the Hilbert-space weak stability~saga.

\vskip5pt
Take a unitary operator ${U\kern-1pt\in\BH}$ on a Hilbert space $\H.$ Since
$U$ is hyponormal (and also a contraction), ${\N(I-U)}$ reduces it, and hence
$U$ is decomposed on $\H={\N(I-U)\oplus\N(I-U)^\perp}\!$ as $U={I\oplus W}\!$,
where $I=U|_{\N(I-U)}$ is the identity on ${\N(I-U)}$ and
$W\!=U|_{\N(I-U)^\perp\!}$ is unitary on ${\N(I-U)^\perp\!}.$ Thus the next
result from \cite[Corollary 3.8]{CJJS2} becomes a particular case of
Corollary~5.10.

%%%%%%%%%%%%%%%%%%%%%%%%%%  PROPOSITION 9.1
\vskip5pt\noi
{\bf Proposition 9.1} \cite{CJJS2} (2024).\
{\it If\/ $U$ is a unitary operator, then
$$
U=I\oplus W,
$$
\vskip-2pt\noi
on\/ ${\H=\N(I-U)\oplus\N(I-U)^\perp}\!$, where\/ $I$ is the identity on\/
${\N(I-U)}$ and\/ $W\!$ is~unitary on\/ ${\N(I-U)^\perp\!}.$ Moreover,
$$
U^n\!\wconv P
\quad\hbox{\it if and only if}\quad
W^n\!\wconv O,
$$
so that the weak limit\/ ${P\in\BH}$ of\/ $\{U^n\}$, if it exists, is the
orthogonal projection
$$
P=I\oplus O
$$
on\/ ${\H=\N(I-U)\oplus\N(I-U)^\perp}\!$ with}\/ $\R(P)=\N({I-U})$.

\vskip5pt
We say that a unitary operator is {\it completely nonidentity}\/ if the
restriction of it to every reducing subspace is not an identity operator$.$
So Proposition 9.1 says~that
\vskip-3pt\noi
$$
\vbox{\noi
{\narrower\narrower
{\it the power sequence of a unitary operator is weakly convergent \\
if and only if its completely nonidentity part is weakly stable}\/.
\vskip0pt}
\vskip-20pt
}                                                                \eqno{(9.1)}
$$
\vskip10pt
\vskip0pt\noi
Note that both power sequences $\{U^n\}$ and $\{W^n\}$ may not converge; for
instance, let $U$ be a symmetry (i.e., a unitary involution) given by
${U=\big(\smallmatrix{1 &  0 \cr
                      0 & -1 \cr}\big)}={1\oplus(-1)}.$

\vskip5pt
Let $\A_\TT$ stand for the $\sigma$-algebra of Borel subsets of the unit
circle $\TT$.

%%%%%%%%%%%%%%%%%%%%%%%%%%%  REMARK 9.2
\vskip5pt\noi
{\bf Remark 9.2.}\
Let ${\lambda\!:\A_\TT\to\RR^{\kern-.5pt+}}\kern-1.5pt$ be a finite
positive-real-valued measure on $\A_\TT$, and let $\mu$ be either another
finite positive-real-valued measure on $\A_\TT$ or a finite
positive-$\BH$-valued spectral measure on $\A_\TT.\kern-1pt$ The following
notions will be required next.
\vskip5pt\noi
(a)
Recall that a spectral measure on $\A_\TT$ is a mapping ${E\!:\A_\TT\to\BH}$,
where $\H$ is a complex Hilbert space, such that (i) $E(\Delta)$ is an
orthogonal projection for \hbox{every} ${\Delta\in\A_\TT}$, (ii)
${E(\void)=O}$, the null operator, and ${E(\TT)=I}$, the identity operator,
(iii) ${E(\Delta_1\cap\Delta_2)}={E(\Delta_1)E(\Delta_2)}$ for every
$\kern-1pt{\Delta_1,\Delta_2\in\kern-.5pt\A_\TT}$, and (iv) $E$ is countably
additive (i.e., ${E(\bigcup_k\Delta_k)\kern-1pt=\!\sum_kE(\Delta_k)}$ whenever
$\{\Delta_k\}$ is a countable collection of~pair\-wise disjoint sets in
$\A_\TT).$ Thus ${O\le E(\Delta_1)\le E(\Delta_2)\le I}$ whenever
${\Delta_1\sse\Delta_2}$ in $\A_\TT.$ For notational simplicity, we will use,
from now on, the same symbol for the number zero and the null operator, as
well as for the number one and the identity operator$.$ So a spectral measure
$E$ on $\A_\TT$ is a finite positive measure in~the~sense~that
$$
\hbox{${0=E(\void)\le E(\Delta)\le E(\TT)=1}\;$ for every
$\;{\Delta\kern-.5pt\in\kern-1pt\A_\TT}.$}
$$

\vskip1pt\noi
(b)
Regarding $\lambda$ as a reference measure for $\mu$, also recall that
\begin{description}
\item{$\kern-9pt$(i)$\kern2pt$}
$\mu$ is absolutely continuous with respect to $\lambda$ if all sets of
$\lambda$-measure$\phantom{.^{\textstyle|}}\kern-3pt$zero have $\mu$-measure
zero, that is, for ${\Delta\in\A_\TT}$,
$$
\lambda(\Delta)=0
\quad\hbox{implies}\quad
\mu(\Delta)=0;
$$
\item{$\kern-11pt$(ii)$\kern1.5pt$}
$\mu$ and $\lambda$ are equivalent if they are absolutely continuous with
respect to each other, that is, for ${\Delta\in\A_\TT}$,
$$
\lambda(\Delta)=0
\quad\hbox{if and only if}\quad
\mu(\Delta)=0;
$$
\item{$\kern-12pt$(iii)$\kern0pt$}
$\mu$ and $\lambda$ are (mutually) singular if they have disjoint supports
or, equivalently, if $\mu$ is concentrated on a set of $\lambda$-measure zero
(and vice versa), which means that there exists a partition
$\{\Delta,\Lambda\}$ of $\TT$ with ${\Delta,\Lambda\in\A_\TT}$ such that
$$
\mu(\Delta)=\lambda(\Lambda)=0;
$$
\item{$\kern-12pt$(iv)$\kern-0pt$}
$\mu$ is continuous with respect to $\lambda$ if for every singleton
${\{\alpha\}\in\A_\TT}$,
$$
\lambda(\{\alpha\})=0
\quad\hbox{implies}\quad
\mu(\{\alpha\})=0;
$$
\item{$\kern-11pt$(v)$\kern1pt$}
$\mu$ is discrete with respect to $\lambda$ if $\mu$ is concentrated
on a countable set of $\lambda$-measure zero, which means that there exists a
partition $\{\Delta,\Lambda\}$ of $\TT$ with ${\Delta,\Lambda\in\A_\TT}$ such
that $\Lambda$ is countable (and all subsets of $\Lambda$ lie in $\A_\TT$) and
$$
\mu(\Delta)=\lambda(\Lambda)=0.
$$
\end{description}
Clearly, if $\mu$ is absolutely continuous, or discrete, with respect to
$\lambda$, then it is continuous, or singular, respectively, with respect to
$\lambda$.
\vskip5pt\noi
(c)
Also recall that if $\lambda$ and $\mu$ are finite measures (actually,
$\sigma$-finiteness is~enough), then $\mu$ has a unique decomposition
${\mu=\mu_a\kern-1pt+\mu_s}$ and a unique decomposition
$\mu={\mu_c\kern-1pt+\mu_d}$, both with respect to a reference measure
$\lambda$, where the finite measures $\mu_a$, $\mu_s$, $\mu_c$, and $\mu_d$
are absolutely continuous, singular, continuous, and discrete, respectively,
with respect to $\lambda.$ So $\mu_s={\mu_{sc}\kern-1pt+\mu_{sd}}$, where
$\mu_{sc}$ and $\mu_{sd}$ are singular-continuous and singular-discrete (i.e.,
discrete), respectively, with respect to $\lambda.$ Hence (this is the
Lebesgue Decomposition Theorem --- see, e.g.,
\cite[Section~7.3]{EMT2}$\kern.5pt$)
\vskip1pt\noi
$$
\mu={\mu_a\kern-1pt+\mu_s}=\mu_a\kern-1pt+\mu_{sc}\kern-1pt+\mu_{sd}.
$$
\vskip-2pt

\vskip5pt
A unitary operator $U\!$ on a complex Hilbert space $\H$ is {\it absolutely
continuous}, {\it continuous}, {\it singular}, {\it singular-con\-tinuous}, or
{\it singular-discrete} (i.e., {\it discrete}\/) if its spectral measure $E$
is absolutely continuous, continuous, singular, singular-continuous, or
singular-discrete with respect to the normalised Lebesgue measure $\lambda$
on $\A_\TT$ (i.e., $\lambda(\TT)=1$ --- recall:\ ${\sigma(U)\sse\TT}$ is
the support of $E).$ By the Spectral Theorem, \hbox{every} unitary operator
$U\!$ on a Hilbert space $\H$ is uniquely decomposed as the \hbox{direct sums}
\vskip4pt\noi
$$
U\kern-1pt=U_{\!a}\oplus U_{\!s}=\kern1ptU_{\!a}\oplus U_{\!sc}\oplus U_{\!sd}
$$
of an absolutely continuous unitary $\kern-1ptU_{\!a}\kern-1pt$ on
$\kern-1pt\H_a$, a singular unitary $\kern-1ptU_{\!s}\kern-1pt$ on
$\kern-1pt\H_s$, a singular-continuous unitary $\kern-1ptU_{\!sc}\kern-1pt$
on $\kern-1pt\H_{sc}$, and a discrete unitary $\kern-1ptU_{\!sd}\kern-1pt$ on
$\kern-1pt\H_{sd}$ (where any part may be missing), with $\kern-1pt\H$
decomposed into orthogonal direct sums of Hilbert spaces:
\vskip4pt\noi
$$
\H\kern-1pt={\H_a\oplus\H_s}\!={\H_a\oplus\H_{sc}\oplus\H_{sd}}.
$$
The direct sum ${U_{\!a}\oplus U_{\!sc}}$ on ${\H_a\oplus\H_{sc}},$
consisting of the absolutely continuous and the singular-continuous parts of
$U\kern-1pt$, is referred to as the {\it continuous}\/ part of $U\kern-1pt$,
although the nomenclature may be evidently tricky --- every unitary operator
is~con\-tinuous in the sense of being linear and bounded.

\vskip5pt
{\it Warning}\/.\
A scalar spectral measure for a normal operator is a finite
positive-real-valued measure equivalent to the spectral measure of its
spectral decomposition$.$ In particular, a scalar spectral measure for a
unitary operator $U\!$ on a Hilbert space $\H$ is any finite
positive-real-valued measure $\mu$ on $\A_\TT$ equivalent to the spectral
measure $E$ of $U\!$ (i.e., $\mu$ and $E$ are mutually absolutely
continuous)$.\kern-1pt$ The definition of scalar spectral measure does not
depend on separability, but its existence does$:$ if $\H$ is separable, then
there exists a scalar spectral measure for $U\!.$ Indeed, for each ${x\in\H}$,
set ${\mu_{x,x}(\Delta)=\<E(\Delta)x\,;x\>=\|E(\Delta)x\|^2}$ for every
${\Delta\in\A_\TT}.\kern-1pt$ This defines a fi\-nite positive-real-valued
measure $\mu_{x.x}$ on $\A_\TT$ for each ${x\in\H}.$ As is well known, if $\H$
is separable, then there is an ${e\in\H}$ such that ${E(\Delta)e\ne0}$
whenever ${E(\Delta)\ne0}$ (i.e., if ${e\in\N(E(\Delta))}$ for some
${\Delta\in\A_\TT}$, then ${\N(E(\Delta))=\H}$) and, for such an ${e\in\H}$,
the measure $\mu_{e,e}$ is a scalar spectral measure for $U$ (i.e.,
$\mu_{e,e}$ and $E$ are equivalent). See, e.g.,
\cite[Definitions 4.4, 4.5, Lemmas 4.6, 4.7, and Remarks on p.99]{ST2}$.$
The~most common definitions of $U_{\!a}$, $U_{\!s}$, $U_{\!sc}$, and
$U_{\!sd}$ are in terms of a scalar \hbox{spectral} measure, which requires
separability for $\H.$ In order to avoid the assumption of sep\-arability, we
have defined $U_{\!a}$, $U_{\!s}$, $U_{\!sc}$, and $U_{\!sd}$ in terms of the
spectral measure $E$, which boils down to the same thing as defining them in
terms of the scalar spectral measure $\mu_{e,e}$ if $\H$ is separable$.$
(The same approach was \hbox{applied in \cite[Section 5]{JJKS}.)}

\vskip5pt
Along this line, the following results on weak stability of each of the above
parts of a unitary operator were summarised in \cite[Section 3]{Kub3}.

%%%%%%%%%%%%%%%%%%%%%%%%%%%  REMARK 9.3
\vskip5pt\noi
{\bf Remark 9.3.}\
Consider the decomposition $U\!={U_{\!a}\oplus U_{\!sc}\oplus U_{\!sd}}$ of a
unitary operator\/ $U$ on a complex Hilbert space
$\H=\H_a\oplus\H_{sc}\oplus\H_{sd}$.
\begin{description}
\item{$\kern-9pt$\rm(a)}
{\it Absolutely continuous unitaries are always weakly stable}\/
$\;$(i.e., $U_{\!a}^n\!\wconv O$).
\vskip3pt
\item{$\kern-9pt$\rm(b)}
{\it Singular-discrete unitaries are never weakly stable}\/
$\;$(i.e., $U_{\!sd}^n\!\notwconv O$).
\vskip3pt
\item{$\kern-9pt$\rm(c)}
{\it If\/ $U$ is weakly stable, then\/ $U$ is continuous}\/
$\;$(i.e., ${U^n\!\wconv O}\limply{U\!=U_{\!a}\oplus U_{\!sc}}$).
\end{description}
\vskip-2pt

\vskip5pt
This gives a first characterisation of weakly stable unitary operators.

%%%%%%%%%%%%%%%%%%%%%%%%%%%  COROLLARY 9.4
\vskip5pt\noi
{\bf Corollary 9.4.}\
{\it A unitary operator is weakly stable if and only if its
singular-contin\-uous part is weakly stable\/ $($if it exists\/$)$ and its
singular-discrete part does not \hbox{exist}}\/.
 
\vskip5pt
Let $\mu$ be a finite positive-real-valued measure on $\A_\TT$, and consider
the Hilbert space $L^2(\TT,\mu)$ of square integrable complex functions on
$\TT$ with respect to $\mu$, so~that
$$
L^2(\TT,\mu)
\;\;\hbox{is separable}.
$$
In fact, as $\TT\kern-1pt$ {\it is compact}\/, the Stone--Weierstrass Theorem
ensures that the set$P(\TT)$of polynomials
${q(\,\cdot\,,\cdot\,)\!:
\kern-1.5pt\TT\kern-2pt\times\kern-2pt\TT\kern-2pt\to\CC}$
in $z$ and $\overline z$ is dense in the space $C(\TT)$ of continuous
functions on $\TT\kern-1pt$ equipped with the sup-norm ${\|\cdot\|_\infty}.$
So $P(\TT)$ is dense in ${(C(\TT),\|\cdot\|_p)}$ for \hbox{every} ${p\ge1}$
because $\mu$ {\it is finite}\/$.$ However, $C(\TT)$ is dense in
${(L^p(\TT,\mu),\|\cdot\|_p)}$ (as Borel measures are regular on compact
metric spaces; see, e.g., \cite[Theorem 29.14]{Bau}$\kern.5pt$).\ Hence
$P(\TT)$ is dense in ${(L^p(\TT,\mu),\|\cdot\|_p)}.$ As $P(\TT)$ is a
linear span of a countable set, ${(L^p(\TT,\mu),\|\cdot\|_p)}$ is
spanned by a countable set, thus separable.\

\vskip5pt
Now consider the unitary multiplication operator $U_{\vphi,\mu}\!$
on $L^2({\TT,\mu})$ induced by the identity function ${\vphi\!:\TT\!\to\!\TT}$
(\hbox{${\vphi(z)=z}$ for ${z\in\TT}$), given by}
$$
(U_{\vphi,\mu}\,\psi)(z)=\vphi(z)\psi(z)=z\psi(z)
\quad\hbox{$\mu$-a.e., for $z\in\TT$},
\quad\hbox{for every $\psi\in L^2(\TT,\mu)$}.
$$
It is well known that the measure $\mu$ is identified with the scalar spectral
measure of the unitary multiplication operator
${U_{\vphi,\mu}\!\in\B[L^2(\TT,\mu)]}$ (see, e.g.,
\cite[\hbox{Remark~p.99}]{ST2}$\kern.5pt$).

\vskip5pt
A finite measure ${\mu\!:\A_\TT\!\to\RR^{\kern-.5pt+}}\kern-1.5pt$ is a
{\it Rajchman measure}\/ if
$$
\hbox{${\int_\TT\!z^k\,d\mu\to0}\;$ as $\;{|k|\to\infty}.$}
$$
With $U_{\vphi,\mu}\kern-1pt$ on ${L^2(\TT\!,\mu)}$ being the unitary
multiplication operator defined~above,
$$
\hbox{\it $\mu$ is a Rajchman measure}
\;\;\;\hbox{\it if and only if}\;\;\;
{U_{\vphi,\mu}}^{\kern-12pt n\kern8pt}\!\!\wconv O.              \eqno{(9.2)}
$$
Indeed, if $\mu$ is Rajchman, then
${\<U^n_{\vphi,\mu}z^k;z^\ell\>=\int_\TT\!z^{n+k-\ell}\,d\mu\to 0}$ as
${n\to\infty}$ for each ${k,\ell}.$ Therefore, since the Stone--Weierstrass
Theorem ensures that the set of polynomials is dense in
${(L^1(\TT,\mu),\|\cdot\|_1)}$, this implies
${\<U^n_{\vphi,\mu}\psi;\psi\>}={\!\int_\TT z^n|\psi|^2d\nu\to0}$ for every
${\psi}$ in ${L^2(\TT,\mu)}$, which in turn implies
${\<U^n_{\vphi,\mu}\psi;\phi\>\to0}$ for \hbox{every} ${\psi,\phi}$ in
${L^2(\TT,\mu)}$ because the Hilbert space ${L^2(\TT,\mu)}$ is complex.\ So
if $\mu$ is Rajchman, then ${U^n_{\vphi,\mu}\!\wconv O}.$ The converse is
straightforward$:$ if\/ ${U^n_{\vphi,\mu}\!\!\wconv O}$, then
${\int_\TT\!z^n\,d\mu(z)}
\kern-1pt=\kern-1pt{\<{U_{\vphi,\mu}}^{\kern-12pt n\kern8pt}1\,;1\>\to0}$
as ${n\kern-1pt\to\kern-1pt\infty}$, and hence
${\int_\TT\!z^k\,d\mu(z)\to0}$ as~${|k|\to\infty}$.\ (The equivalence in (9.2)
was applied in \cite[p$.$1384, proof of Theorem 5]{BM},
\cite[Proposition 3.3]{Kub3}, and \cite[Eq$.$~(5.3)]{JJKS}.)\

\vskip5pt
Properties of Rajchman measures that will be required here are listed below
(see, e.g., \cite[p.364 and \hbox{Theorem}~3.4]{Lyo}$\kern.5pt).$ The terms
`absolutely continuous'$\kern-1pt$, `\hbox{continuous}'$\kern-1pt$,
`singular'$\kern-1pt$, and `discrete' are with respect to the normalised
Lebesgue \hbox{measure on $\kern-1pt\A_\TT\kern-1pt$}.
\vskip5pt\noi
$$
\vbox
{\noi
{\narrower\narrower
{\it Every absolutely continuous finite measure is Rajchman}\/. \\
{\it Every Rajchman measure is continuous}\/. \\
{\it There exist singular Rajchman measures}\/.
\vskip3pt\noi
Therefore, {\it every singular Rajchman measure is singular-continuous}\/.
\vskip0pt}
\vskip-34pt
}                                                                \eqno{(9.3)}
$$
\vskip20pt\noi

\vskip5pt
The equivalence in (9.2) says that such a special unitary multiplication
operator $U_{\vphi,\mu}\!$ is weakly stable if and only if $\mu$ is a Rajchman
measure.\ In light of this, the next results have been verified in
\cite[Corollary 5.3 and Theorem 7.3]{Kub5}.

%%%%%%%%%%%%%%%%%%%%%%%%%%  PROPOSITION 9.5
\vskip5pt\noi
{\bf Proposition 9.5} \cite{Kub5} (2026).\
{\it If a finite measure\/ $\mu$ is continuous, then\/ $U_{\vphi,\mu}\!$ is
weakly quasistable$.$ Moreover, $U_{\vphi,\mu}\!$ is weakly stable if and
only if}\/ ${\lim_n\kern-1pt\<U^n_{\vphi,\mu}1\,;1\>\kern-1pt=\kern-1pt0}.$

\vskip5pt
According to Corollary 9.4, weak stability of unitary operators is reduced to
weak stability of singular-continuous unitary operators$.$ Examples of weakly
stable and weakly unstable singular-continuous unitary operators are given
below, where the examples in (a) and (b$_1$) have been pulled together in
\cite[Propositions 3.2~and~3.3]{Kub3}.

%%%%%%%%%%%%%%%%%%%%%%%%%%%  EXAMPLE 9.6
\vskip5pt\noi
{\bf Example 9.6}.\
There are stable and unstable singular-continuous unitary operators.
\begin{description}
\item{$\kern-9pt$\rm(a)}
A weakly stable singular-continuous unitary operator.
\vskip4pt\noi
If a finite measure $\mu$ on $\A_\TT$ is a singular Rajchman measure, then
according to (9.3), the multiplication operator $U_{\vphi,\mu}$ appearing in
(9.2) is a weakly stable singular-continuous unitary operator$.$
\vskip4pt
\item{$\kern-14pt$\rm(b$_1$)}
A weakly unstable singular-continuous unitary operator.
\vskip4pt\noi
If $\lambda$ is the normalised Lebesgue measure on $\A_\TT$, then $U\!$ on
${L^2(\TT,\lambda)}$ given by
\vskip4pt
\centerline{
$(U\!\psi)(z)=z^q\psi(\gamma z)$
\quad $\lambda$-a.e., for $z\in\TT$,
\quad for every $\psi\in L^2(\TT,\lambda)$}
\vskip4pt\noi
is a unitary operator$.$ If $q$ is a sufficiently small nonzero rational
(e.g., $0<|q|\le 1/12$) and $\gamma$ is an irrational in $\TT$ (i.e.,
${\gamma=e^{2\pi i\alpha}}$ with $\alpha$ being an irrational in ${(0,1)}$),
then $U\!$ is weakly unstable (i.e., ${U^n\!\notwconv O}$) and
singular-continuous$.$
\vskip4pt
This example appeared in \cite[Answer 3]{Kub1}$.$ It was based on two facts:\
(i) $U\!$ is not discrete and (ii) there is at least one subsequence
$\{U^{n_k}\}$ of $\{U^n\}$ for which ${0<\inf_k|\<U^{n_k}1\,;1\>|}$
\cite[Lemma]{Cho}$.$ Then $U$ is not weakly stable, and so not absolutely
continuous$.$ Also, its scalar spectral measure is pure:\ either purely
absolutely continuous, purely singular-continuous, or purely discrete
\cite[\hbox{Theorem 3}]{Hel}; so it is singular-continuous.\
(See \cite{Cho2}~as~well.)
\vskip4pt\noi
\item{$\kern-14pt$\rm(b$_2$)}
Another weakly unstable singular-continuous unitary operator.
\vskip4pt\noi
A Borel--Stieltjes measure $\mu_{[0,1]}$ on the $\sigma$-algebra $\A_{[0,1]}$
of Borel subsets of ${[0,1]}$, generated by the Cantor function associated
with the Cantor set over ${[0,1]}$, is called the Cantor--Lebesgue measure
over ${[0,1]}$.\ The Cantor--Lebesgue measure $\mu$ over $\TT$ is the measure
on $\A_\TT$ induced by Cantor--Lebesgue meas\-ure $\mu_{[0,1]}$ over ${[0,1]}$
by the function ${\gamma\!:[0,1]\to\TT}$, given by
${\gamma(\alpha)=e^{2\pi i\alpha}}$ for ${\alpha\in[0,1)}$, such that
${\mu(\Delta)=\mu_{[0,1]}(\gamma^{-1}(\Delta))}$ for ${\Delta\in\A_\TT}.$
It is known~that
\vskip4pt
\centerline{the Cantor--Lebesgue measure on $\A_{[0.1]}$ is
singular-continuous}
\vskip4pt\noi
(with respect to normalised Lebesgue measure on $\A_{[0.1]}$; see, e.g.,
\cite[Problem 7.15(c)]{EMT2} or \cite[Example 3, Section 1.4]{RS}$\kern.5pt$),
and so is the induced $\mu$~by~$\gamma$:
\vskip4pt
\centerline{the Cantor--Lebesgue measure on $\A_\TT$ is singular-continuous}
\vskip4pt\noi
(with respect to normalised Lebesgue measure on $\A_\TT$, the so-called
arc-length measure on $\A_\TT$).\ Moreover, it is also known that (see, e.g.,
\cite[p.364]{Lyo}$\kern.5pt$)
\vskip4pt
\centerline{the Cantor--Lebesgue measure on $\A_\TT$ is not a Rajchman
measure.}
\vskip4pt\noi
So we get the following singular-continuous weakly unstable unitary
\hbox{operator}:
\vskip4pt\noi
{\narrower
if $\mu$ is the Cantor--Lebesgue measure on $\A_\TT$, then the multiplication
$U_{\vphi,\mu}$ in (9.2) is a weakly unstable singular-continuous unitary
operator.
\vskip0pt}
\end{description}
\vskip-3pt

\vskip5pt
It was asked in \cite{KD} if there is a singular-continuous unitary operator
that is not weakly quasistable.\ In particular, are those weakly unstable
unitary operators in Examples 9.6 (b$_1$ and b$_2$) weakly quasistable\/$?$
(See also \cite[Proposition IV.1.2(b)]{Eis} for a result along this line
regarding almost weak stability.)

\vskip5pt
Let ${E\!:\A_\TT\!\to\BH}$ be the spectral measure of a unitary operator
$U\!$ on a \hbox{complex} Hilbert space $\H.$ For each ${x,y\in\H}$
take the complex measure ${\mu_{x,y}\!:\A_\TT\to\CC}$ given by
$\mu_{x,y}(\Lambda)={\<E(\Lambda)x\,;y\>}$ for every
${\Lambda\in\kern-1pt\A_\TT}.$ By the Spectral Theorem,
${U\!=\!\int_\TT\kern-1pt z\,dE_\lambda}$, and so
${U^n\!=\!\int_\TT\kern-1pt z^ndE_\lambda}$, which means
${\<U^nx\,;y\>}=\!{\int_\TT\kern-1pt z^n d\<E_\lambda x\,;y\>}.$ Hence
${\<U^nx\,;x\>}=$ $\!{\int_\TT\kern-1pt z^n d\<E_\lambda x\,;x\>
=\kern-1pt\int_\TT\kern-1pt z^nd\mu_{x,x}}.$
As $\H$ is complex, this establish the next result:
$$
\hbox{${U^n\!\wconv O}\;$ if an only if
$\;{\int_\TT\kern-1pt z^nd\mu_{x,x}\!\to0}\;$ for every $\;{x\in\H}$.}
$$
\vskip-2pt

%%%%%%%%%%%%%%%%%%%%%%%%%%%  PROPOSITION 9.7
\vskip5pt\noi
{\bf Proposition 9.7} \cite{JJKS} (2024).\
{\it $\!$A unitary operator $U\kern-3pt$ on $\H\!$ with spectral measure $E$
is weakly stable if and only if $\kern.5pt{\<E(\cdot)x\,;x\>}$ is a
Rajchman measure for every}\/ ${x\in\H}.$

\vskip5pt
Proposition 9.7 extends (9.2), and it can also be viewed as consequence of
(9.2)$.$ In fact, the essence of the Spectral Theorem is that any normal
operator is unitarily equivalent to a (possibly uncountable) orthogonal direct
sum of cyclic multiplication normal operators (on separable Hilbert spaces),
and a sum (possibly uncountable) of Rajchman measures is again a Rajchman
measure$.$ Proposition 9.7 was originally proved in \cite[Corollary 6.2]{JJKS}
as a particular case of Proposition 5.13(a).\ (See also
\cite[Section IV.1]{Eis} and the references therein for Rajchman measures
applied to~stability.)

\vskip5pt
The result below is from \cite[Lemmas 5.3, 5.4, Theorem 5.5]{JJKS}$.$ It
extends the weak stability criterion of Corollary 9.4, and is a further
consequence of \hbox{Corollary 5.10}.

%%%%%%%%%%%%%%%%%%%%%%%%%%%  PROPOSITION 9.8
\vskip5pt\noi
{\bf Proposition 9.8} \cite{JJKS} (2024).\
{\it Let\/ $U\!$, $U_a$, $U_{sc}$, and $U_{sd}$ be unitary operators, where
$U_a$ is absolutely continuous, $U_{sc}$ is singular-continuous, and\/
$U_{sd}$ is singular-discrete}\/.
\begin{description}
\item{$\kern-9pt$\rm(a)}
{\it The power sequence of a singular-continuous unitary operator converges
weakly if and only if it is weakly stable$.$ That is}\/,
\goodbreak\vskip4pt\noi
$$
{U_{sc}}^{\kern-6pt n\kern2pt}\!\wconv P
\quad\hbox{\it if and only if}\quad
{U_{sc}}^{\kern-6pt n\kern2pt}\!\wconv O.
$$
\item{$\kern-9pt$\rm(b)}
{\it The power sequence of a singular-discrete unitary operator converges
weakly if and only if it is the identity operator$.$ That is}\/,
$$
{U_{sd}}^{\kern-6pt n\kern2pt}\!\wconv P
\quad\hbox{\it if and only if}\quad
U_{sd}=I.
$$
\item{$\kern-9pt$\rm(c)}
{\it Regarding the decomposition\/
$U\kern-3pt=\kern-1pt{U_a\kern-1.5pt\oplus U_{sc}\kern-1.5pt\oplus U_{sd}}$
for a unitary operator\/ $U\!$~on~$\H$,}
$$
U^n\!\wconv P
\quad\hbox{\it if and only if}\quad
{U_{sc}}^{\kern-6ptn}\wconv O,
\quad{U_{sd}}=I,
$$
{\it or any of the parts\/ $U_{sc}$ or\/ $U_{sd}$ are absent in the above
decomposition, where\/ $P$ is the orthogonal projection on the Hilbert space
${\H\kern-1pt=\kern-1pt
\H_a\kern-1pt\oplus\kern-1pt\H_{sc}\kern-1pt\oplus\kern-1pt\H_{sd}}$
given~by
$$
P=O\oplus O\oplus I,
$$
where, again, any part in the above direct sum may be absent}\/.
\end{description}
\vskip-1pt

\vskip5pt
Reducing subspaces were defined in Section 2, which is a property that makes
sense only in a Hilbert-space setting$.$ A Hilbert-space operator is (i)
reducible if it has a nontrivial reducing subspace; otherwise, it is
irreducible; and it is (ii) reductive if every invariant subspace reduces it;
otherwise, it is nonreductive$.$ Examples$:$ \hbox{every} self-adjoint
operator is trivially reductive ($\M$ reduces $T$ if and only if it is
invariant for both $T$ and $T^*$), and every compact normal operator is also
reductive \cite[Theo\-rem 1]{And}, but there are nonreductive unitary
operators (e.g., bilateral shifts$:$~every uni\-lateral shift is the
restriction of a bilateral shift to an invariant subspace)$.$ The question
that asks whether there is a reductive operator that is not normal was shown
in \cite[Corollary 3.1]{DPP} to have an affirmative answer if and only if
every operator on a Hilbert space has a nontrivial invariant subspace.\

\vskip5pt
We close the paper by considering the relationship among the shape of the
spectrum, reductivity, and weak stability for unitary operators.\ Such an
investigation has its starting point in the following facts:
\begin{description}
\item{$\kern-12pt${\sc(i)}$\kern2pt$}
{\it a unitary operator is absolutely continuous if and only if it is a part
of a bilateral shift}\/ (or a bilateral shift itself)
\cite[Exercise 6.8, p.56]{Fil}$\kern.5pt$),
\vskip4pt
\item{$\kern-14pt${\sc(ii)}$\kern2pt$}
{\it a unitary operator is reductive if and only if no part of it is a
bilateral shift}\/
(see, e.g., \cite[1.VI, p.18]{Fil} or \cite[Proposition 1.11]{RR}$\kern.5pt$),
\vskip4pt
\item{$\kern-15.5pt${\sc(iii)}$\kern0pt$}
{\it a singular unitary operator is reductive}\/ (by the above two results),
\vskip4pt
\item{$\kern-14pt${\sc(iv)}$\kern1pt$}
{\it if the spectrum of a unitary operator is not the whole unit circle, then
it is reductive}\/
(see, e.g., \cite[Theorem 13.2]{Fil} or \cite[Theorem 1.23]{RR}, and the~
converse fails --- see Example 9.9(c) below).
\end{description}
\vskip-2pt

\vskip5pt\noi
Corollary 9.4 and Example 9.6 naturally lead us to enquire:
\vskip5pt\noi
\centerline{\it when is a singular-continuous unitary operator weakly
stable\,$?$}
\vskip5pt\noi
Recall that the spectrum of a unitary operator is included in the unit
circle.\ A possible answer to the above question does not depend on whether or
not the spec\-trum is the whole unit circle or whether the unitary operator
is reductive.

%%%%%%%%%%%%%%%%%%%%%%%%%%%  EXAMPLE 9.9
\vskip5pt\noi
{\bf Example 9.9.}\
Consider the possibilities$:$ a unitary operator is weakly stable or
not, and its spectrum is the whole unit circle or not$.$ All combinations
\hbox{are possible}.
\begin{description}
\item{$\kern-9pt$\rm(a)}
As we saw above, {\it a bilateral shift\/ $S$ on a Hilbert space is a
classical example of an absolutely continuous $($nonreductive\/$)$ unitary
operator}\/, {\it thus weakly stable}\/ (${S^n\!\wconv O}$ as in
(4.5)$\kern.5pt$), {\it whose spectrum is the unit circle}\/
(${\sigma(S)=\TT})$.
\item{$\kern-9pt$\rm(b)}
{\it There are absolutely continuous unitary operators whose spectra
are$\phantom{.^{\textstyle|}}\kern-3pt$not
the whole unit circle, which are weakly stable}\/$.$ For instance, let
$\lambda$ be the normalised Lebesgue measure on $\A_\TT$, and take the
separable Hilbert space $L^2(\TT,\lambda)$ for which $\{e_k\}_{k\in\ZZ}$, with
${e_k(z)=z^k}$ for\/ ${z\in\TT}$, is an orthonormal basis$.$ Let
\hbox{$U_{\vphi,\lambda}\!$ on} $L^2({\TT,\lambda})$ be the multiplication
operator induced by the identity function $\vphi.$ So $U_{\vphi,\lambda}$ is
a bilateral shift on $L^2(\TT,\lambda)$, shifting the orthonormal basis
$\{e_k\}_{k\in\ZZ}.$ Thus it is unitary and
${\sigma(U_{\vphi,\lambda})\kern-1pt=\kern-1pt\TT}$, weakly stable
(${{U_{\vphi,\lambda}}^{\kern-12pt n\kern8pt}\!\!\wconv O})$ as in (4.5), and
absolutely continuous by {\sc(i)}$.$ If
${\varUpsilon\kern-1pt\in\kern-1pt\A_\TT}$ is a compact set with
${0\kern-1pt<\kern-1pt\lambda(\varUpsilon)\kern-1pt<\kern-1pt1}$, then
$L^2(\varUpsilon,\lambda)$ is identified (i.e., unitarily equivalent) to a
subspace of $L^2(\TT,\lambda)$, which reduces the absolutely continuous
unitary $U_{\vphi,\lambda}$ (see, e.g., \cite[Theorem 3.6]{RR}$\kern.5pt)$,
and so $U\!={U_{\vphi,\lambda}|_{L^2(\varUpsilon,\,\lambda)}}$ is an
absolutely continuous unitary by {\sc(i)}, with
${\sigma(U)=\varUpsilon\kern-.5pt\subset\kern-1pt\TT}$ (proper inclusion, so
$U$ is {\it reductive}\/ by {\sc(iv)}), and weakly stable (${U^n\!\wconv O})$
by Remark 4.1 or by Remark 9.3(a).
\vskip4pt
\item{$\kern-9pt$\rm(c)}
{\it There are singular-discrete unitary operators whose spectra are the unit
\hbox{circle}, which are not weakly stable}\/$.$ Indeed, if
$\{\alpha_k\}_{k\in\NN}$ is an enumeration of ${\QQ\cap[0,1)}$, then the
unitary diagonal
${D\kern-1pt=\kern-1pt\diag(\kern-.5pt\{e^{2\pi i\alpha_k}\kern-.5pt\})}$ on
$\ell_+^2\kern-1pt$ is singular-discrete\/$.$ In~fact, $D$ is a
{\it reductive}\/ diagonal \cite[Example 13.5]{Dow} with
${\sigma_{\kern-1pt P}(D)=\kern-1pt\{e^{2\pi i\alpha_k}\}_{k\in\NN}}.$
Actu\-ally, no part of $D$ is a bilateral shift --- cf.\ {\sc(ii)}.\ So all
parts of $D$ are not contin\-uous $($see, e.g.,
\cite[Theorem 12.29]{Rud}$\kern.5pt).$ Thus $D$ is not absolutely~continuous,
and hence $D$ is singular-discrete by the Lebesgue Decomposition~Theorem$.$
Then $D$ is weakly unstable $({D^n\!\notwconv O})$ by Remark 9.3(b).\ Since
$D$ is unitary with $\sigma_{\kern-1ptP}(D)$ dense in $\TT$, its spectrum is
the unit circle $({\sigma(D)=\TT})$.
\vskip4pt
\item{$\kern-9pt$\rm(d)}
{\it There are singular-continuous unitary operators whose spectra are not the
whole unit circle, which are weakly unstable}\/$.$ In fact, let
${U_{\vphi,\mu}\!\in\kern-.5pt L^2[\TT,\mu]}$ be~the unitary multiplication
operator induced by the identity function $\vphi$ on $\TT$,~where $\mu$ is the
Cantor--Lebesgue measure on $\A_\TT\kern-1pt$ which, as defined in
Example~9.6(b$_2$), is the measure on $\A_\TT$ induced by the Cantor--Lebesgue
measure $\mu_{[0,1]}$ on $\A_{[0,1]}$ generated by the Cantor function,
associated with the Cantor set~$C$~in ${[0,1]}.$ So the support of
$\mu_{[0,1]}$ is $C$ itself.\ Thus the support of the
\hbox{Cantor}--\hbox{Lebesgue} measure $\mu$ on $\A_\TT$ is the image
${\varGamma\!=\gamma(C)}$ of $C$ under the function ${\gamma\!:[0,1]\to\TT}$,
and therefore $\varGamma\kern-1pt$ is a proper subset of the unit circle
$\TT\kern-1pt.\kern-1pt$ We saw in Example 9.6(b$_2$) that $U_{\vphi,\mu}$ is
singular-continuous and weakly unstable
$({{U_{\vphi,\mu}}^{\kern-12pt n\kern6pt}\!\notwconv O}).$ Its spectrum is the
support of its scalar spectral measure, which is identified with $\mu$, so
that $\sigma(U_{\vphi,\mu})={\varGamma\subset\TT}$ (proper inclusion, and
hence $U_{\vphi,\mu}$ is {\it reductive}\/).
\end{description}

%%%%%%%%%%%%%%%%%%%%%%%%%%%%%%%%%%%%%%%%%%%%%%%%%%%%%%%%%%  ACKNOWLEDGEMENT
\section*{Acknowledgement}

Of course, this expository survey was not written overnight.\ I had the help
of my coauthors on weak stability over the years, who have been cited in the
reference list below.\ I am really grateful to all of them for everything I
have learnt from~them.

%%%%%%%%%%%%%%%%%%%%%%%%%%%%%%%%%%%%%%%%%%%%%%%%%%%%%%%%%  REFERENCES
\bibliographystyle{amsplain}

\end{document}